\documentclass{amsart}
\usepackage{amsmath,amsthm,bbm}
\usepackage{amsfonts,amssymb}

\newtheorem{thm}{Theorem}[section]
\newtheorem{cor}[thm]{Corollary}
\newtheorem{lem}[thm]{Lemma}
\newtheorem{prop}[thm]{Proposition}

\newtheorem{defn}[thm]{Definition}
\newtheorem{re}[thm]{Remark}

\theoremstyle{remark}

\def\1{{\mathbf 1}}

\def\CD{{\mathcal D}}

\def\CL{{\mathcal L}}
\def\CM{{\mathcal M}}
\def\CO{{\mathcal O}}

\def\CS{{\mathcal S}}

\def\CC{{\mathbb C}}
\def\KK{{\mathbb K}}
\def\NN{{\mathbb N}}
\def\PP{{\mathbb P}}
\def\QQ{{\mathbb Q}}
\def\RR{{\mathbb R}}

\def\supp{\operatorname{supp}}

\newcommand{\wt}{\widetilde}
\newcommand{\wh}{\widehat}

\def\cL{{\mathcal L}}
\def\cM{{\mathcal M}}

\def\cS{{\mathcal S}}

\def\cX{{\mathcal X}}

\newcommand{\eps}{{\varepsilon}}

\def\R{{\mathbb R}}

\def\a{{\alpha}}
\def\b{{\beta}}
\def\g{{\gamma}}
\def\ha{{\wh a}}
\def\ha{{\widehat a}}
\def\hb{{\widehat b}}

\def\PP{{\mathcal H}}
\def\ph{{\varphi}}

\def\cS{{\mathcal S}}

\def\cM{{\mathcal M}}
\def\ONE{{\mathbbm 1}}

\def\L2{{L^2}}

\def\wtd{\widetilde}

\def\V{V}
\def\W{W}
\def\LL{\Lambda}
\def\tLL{\widetilde\Lambda}

\def\one{\bf 1}
\def\PsiP{{\mathcal F}}
\def\cF{{\mathcal F}}

\def\WW{{W_\alpha}}
\def\r{{\rho}}
\def\Fsrpq{{F_{pq}^{s\r}}}
\def\fsrpq{{f_{pq}^{s\r}}}
\def\ONE{{\mathbbm 1}}
\def\tONE{{\tilde{\mathbbm 1}}}

\def\Bsrpq{{B_{pq}^{s\r}}}

\def\bsrpq{{b_{pq}^{s\r}}}

\newcommand\F[4]{F^{#1#2}_{#3#4}}
\newcommand\B[4]{B^{#1#2}_{#3#4}}

\def\w{{w_\alpha}}

\def\one{{\bf 1}}

\def\dm{{\bigstar}}
\def\cd{c_\star}

\def\tm{\widetilde{m}}
\def\diam{{\rm diam}\,}
\def\PP{\cF^\a}
\def\expon{\varrho}
\def\KK{\mathcal K}
\def\ss{s}
\def\DD{\partial}
\def\UU{{\mathcal U}}
\def\QQ{{\mathcal Q}}
\def\AA{{\mathcal A}}

\begin{document}

\title{Decomposition of Triebel-Lizorkin and Besov spaces
in the context of Laguerre expansions}

\author{G. Kerkyacharian, P. Petrushev, D. Picard, and Yuan Xu}

\address{Laboratoire de Probabilit\'{e}s et Mod\`{e}les Al\'{e}atoires, CNRS-UMR 7599,
Universit\'{e} Paris VI et Universit\'{e} Paris VII, rue de Clisson, F-75013 Paris}
\email{kerk@math.jussieu.fr}

\address{Department of Mathematics\\University of South Carolina\\
Columbia, SC 29208\\
and Institute of Mathematics and Informatics, Bulgarian Academy of Sciences}
\email{pencho@math.sc.edu}

\address{Laboratoire de Probabilit\'{e}s et Mod\`{e}les Al\'{e}atoires, CNRS-UMR 7599,
Universit\'{e} Paris VI et Universit\'{e} Paris VII, rue de Clisson, F-75013 Paris}
\email{picard@math.jussieu.fr}

\address{Department of Mathematics\\ University of Oregon\\
    Eugene, Oregon 97403-1222.}\email{yuan@math.uoregon.edu}

\date{April 28, 2008}


\thanks{The second author has been supported by NSF Grant DMS-0709046
and the fourth author by NSF Grant DMS-0604056.}

\begin{abstract}
A pair of dual frames with almost exponentially localized elements (needlets) are constructed
on $\RR_+^d$ based on Laguerre functions.
It is shown that the Triebel-Lizorkin and Besov spaces induced by Laguerre
expansions can be characterized in terms of respective sequence spaces that involve
the needlet coefficients.
\end{abstract}

\maketitle

\pagestyle{myheadings}
\thispagestyle{plain}
\markboth{G. KERKYACHARIAN, P. PETRUSHEV, D. PICARD, AND YUAN XU}
{LAGUERRE TRIEBEL-LIZORKIN AND BESOV SPACES}

\section{Introduction}
\setcounter{equation}{0}

The primary goal of this paper is to construct frames on $\RR_+^d:=(0, \infty)^d$
with nearly exponentially localized elements,
based on Laguerre functions and utilize them to the characterization of spaces of
distribution on $\RR_+^d$.
We are interested in extending the fundamental results of Frazier and Jawerth
\cite{F-J1, F-J2, F-J-W} on the $\varphi$-transform on $\RR^d$
in the context of Laguerre expansions.

From the three types of Laguerre functions available in the literature we focus our
attention on the Laguerre functions $\{\cF^\a_\nu\}$ (see (\ref{def.Psi-n}))
which form an orthonormal basis
for the space $L^2(\RR^d_+, \w)$ with weight
\begin{equation}\label{weight}
\w(x):=\prod_{j=1}^d x_j^{2\alpha_j+1}. 
\end{equation}
For various technical reasons we will assume that $\alpha_j\ge 0$,
while in general $\alpha_j>-1$.
The other two classes of Laguerre functions $\{\cL^\a_\nu\}$ and $\{\cM^\a_\nu\}$
(see (\ref{def.dL_n})-(\ref{def.dM_n})) form orthogonal bases for $L^2(\RR^d_+)$ (weight 1).
The d-dimensional Laguerre functions $\cF^\a_\nu$ are products of
univariate Laguerre functions, namely,
$\cF^\a_\nu(x):= \cF^\a_{\nu_1}(x_1)\cdots \cF^\a_{\nu_d}(x_d)$
(see (\ref{def.Psi-n}), (\ref{def.dPsi_n})).
Hence the kernel of the orthogonal projector onto
\begin{equation}\label{def-Wn}
W_n:={\rm span} \{\cF^\a_\nu: |\nu|=n\}
\quad \mbox{is given by}\quad
\cF^\a_n(x, y):= \sum_{|\nu|=n} \cF^\a_\nu(x)\cF^\a_\nu(y).
\end{equation}
Denote $V_n:= \bigoplus_{m=0}^n W_m$. Evidently,
$K_n(x, y):=\sum_{m=0}^n\cF^\a_m(x, y)$
is the kernel of the orthogonal projector onto $V_n$.
A main point in the present paper is that
for compactly supported $C^\infty$ cut-off functions $\ha$
which are constant around zero the kernels
\begin{equation}\label{defn-Lambda-n}
\LL_n(x,y):= \sum_{j=0}^\infty \ha\Big(\frac{j}{n}\Big) \cF^\a_j(x,y)
\end{equation}
decay rapidly (almost exponentially) away from the main diagonal in $\RR^d_+$
(Theorem~\ref{thm:est-Lambda-n}).
For the same kind of kernels associated with the Laguerre functions $\{\cM^\a_\nu\}$
in dimension $d=1$ this fact is established in \cite{Epp2}.
We show that similar results are valid for $\{\cM^\a_\nu\}$ and $\{\cL^\a_\nu\}$
in dimension $d>1$ as well.

We utilize the kernels from (\ref{defn-Lambda-n}) to the construction of
a pair of dual frames $\{\ph_\xi\}_{\xi\in\cX}$ and $\{\psi_\xi\}_{\xi\in\cX}$
with $\cX$ a multilevel index set.
As in other similar settings, the almost exponential localization of
$\ph_\xi$ and $\psi_\xi$ prompts us to call them ``needlets".
The needlet systems from this paper can be regarded as analogues of
the $\varphi$-transform of Frazier and Jawerth \cite{F-J1, F-J2}.
They are particularly well suited for characterization of the Triebel-Lizorkin
and Besov spaces associated with Laguerre expansions.
To be more precise, let
$\ha\in C^\infty$,
$\supp \ha \subset [1/4, 4]$, and
$|\ha|>c$ on $[1/3, 3]$
and define
$$
\Phi_0(x, y) := \cF_0^\a(x, y)
\quad\mbox{and}\quad
\Phi_j(x, y) := \sum_{m=0}^\infty \ha
\Big(\frac{m}{4^{j-1}}\Big)\cF_m^\a(x,y), ~~~j\ge 1,
$$
Then for all appropriate indices (see Definition~\ref{defn-Tr-Liz})
the Laguerre-Triebel-Lizorkin space $\Fsrpq$ is defined as
the set of all tempered distributions $f$ on $\RR^d_+$ such that
$$
\|f\|_{\Fsrpq}:=\Big\|\Big(\sum_{j=0}^{\infty}
\Big[2^{sj}\WW(4^j;\cdot)^{-\r/d}|\Phi_j*f(\cdot)|\Big]^q\Big)^{1/q}\Big\|_p <\infty.
$$
Here $\Phi_j*f(x):=\langle f, \overline{\Phi_j(x, \cdot)}\rangle$
(Definition~\ref{def:convolution})
and the weight $\WW(n;x)$ is define by
\begin{equation}\label{defn-Wnx}
\WW(n;x):= \prod_{j=1}^d (x_j+n^{-1/2})^{2\alpha_j+1}.
\end{equation}
Just for convenience we use dilations by factors of $4^j$ on the frequency side as opposed
to the traditional binary dilation.
The Laguerre-Besov spaces are defined by the (quasi-)norm
$$
\|f\|_{\Bsrpq} :=
\Big(\sum_{j=0}^\infty \Big(2^{s j}
\|\WW(4^j; \cdot)^{-\r/d}\Phi_j*f(\cdot)\|_p\Big)^q\Big)^{1/q}.
$$
Unlike in the classical case on $\RR^d$ the weight $\w$ creates some
inhomogeneity which compels us to introduce the additional term
$\WW(4^j; \cdot)^{-\rho/d}$ with parameter $\rho\in\RR$.
This allows to consider different scales of Triebel-Lizorkin and Besov spaces.
For instance, a ``classical" choice would be $\rho=0$. However, more natural
to us are the spaces $F_{pq}^{ss}$ and $B_{pq}^{ss}$ which embed ``correctly"
with respect to the smoothness parameter $s$.

The main results in this article assert that the Laguerre Triebel-Lizorkin and Besov spaces
can be characterized in terms of respective sequence spaces involving the needlet
coefficients of the distributions (Theorems~\ref{thm:needlet-Tr-Liz}, \ref{thm:Bnorm-eq}).

Along the same lines one can develop a similar theory on $\RR^d_+$ with weight 1
using the Laguerre functions $\{\cL^\a_\nu\}$ or $\{\cM^\a_\nu\}$.
For such spaces induced by $\{\cL^\a_\nu\}$, see \cite{Dzub}.

This paper is an integral part of a broader undertaking for needlet
characterization of Triebel-Lizorkin and Besov spaces on nonstandard domains
(and with weights) such as the sphere \cite{NPW}, interval \cite{KPX1}, ball \cite{KPX2},
and in the setting of Hermite expansions \cite{PX3}.

The outline of the paper is as follows.
All the information we need about Laguerre polynomials and functions is given in \S2.
The localized kernels induced by Laguerre functions are given in \S3.
Some additional background material is collected in \S4.
The construction of needlets is given in \S5.
In \S6 the Laguerre-Triebel-Lizorkin spaces are introduced and characterized
in terms of needlet coefficients,
while the characterization of the Laguerre-Besov spaces is given in \S7.
Some proofs for Sections 3 - 4 are given in \S8 and for Sections 5 - 6
in \S 9.

The following notation will be used throughout:
$\|x\|:=\max_{i}|x_i|$,
$|x|:=\sum_{i=1}^d |x_i|$,
$\|x\|_2:=\Big(\sum_{i=1}^d |x_i|^2\Big)^{1/2}$,
$\|f\|_p:=\Big(\int_{\RR^d_+}|f(x)|^p\w(x)dx\Big)^{1/p}$;
$|E|$ stands for the Lebesgue measure of $E\subset \RR^d_+$,
$\mu(E):= \int_E\w(x)dx$,
$\ONE_E$ is the characteristic function of $E$, and
$\tONE_E:=\mu(E)^{-1/2}\ONE_E$.
Positive constants are denoted by $c$, $c_1, c_*, \dots$ and they
may vary at every occurrence; $A\sim B$ means $c_1A\le B\le c_2 A$.

\section{Background: Laguerre polynomials and functions}\label{Laguerre-poly}
\setcounter{equation}{0}

In this section we collect the information on Laguerre polynomials and functions
that will be needed in this paper.
The Laguerre polynomials $L_n^\a$ $(\alpha >-1)$
can be defined by their generating function
$$
\sum_{n=0}^\infty L_n^\a (x) r^n = (1-r)^{-\a -1} e^{-xr /(1-r)}, \qquad |r| < 1.
$$
They are orthogonal on $\RR_+ = (0,\infty)$ with weight
$x^\a e^{-x}$ , more precisely,
$$
\int_0^\infty L_n^\a(x)L_m^\a(x)e^{-x}x^\a dx
=\frac{\Gamma(n+\alpha+1)}{\Gamma(n+1)}\delta_{n,m}
=\Gamma(\alpha+1)L_n^\a(0)\delta_{n,m},
$$
where we used that $L_n^\a(0) = \binom{n+\a}{n}$ \cite[(5.1.1)]{Sz}.

Let $L_\nu^\a(x) : =  L_{\nu_1}^{\a_1}(x_1) \cdots  L_{\nu_d}^{\a_d}(x_d)$
be the product Laguerre polynomials on $\RR_+^d$,
where $\nu = (\nu_1,\ldots,\nu_d) \in \NN_0^d$ and $\a = (\a_1,\ldots, \a_d)$.
For $\delta > -1$, define
\begin{equation} \label{cesaro-kernel}
    P_n^{\a,\delta}(x;y) : = \sum_{k=0}^n A_{n-k}^\delta \sum_{|\nu| =k}
           \frac{L_\nu^\a(x)  L_\nu^\a(y)} {L_\nu^\a(0)},
    \qquad A_m^\delta : = \binom{m+\delta}{m}.
\end{equation}
This is a constant multiple of the $n$th Ces\`aro sum of the reproducing kernels for Laguerre
polynomials in dimension $d$.
Using the generating
function of the Laguerre polynomials, it is shown in \cite{X00} that
\begin{equation}\label{cesaro}
     P_n^{\a,\delta}(x , 0) = L_n^{|\a| + \delta + d} (|x|). 
\end{equation}
The product formula for Laguerre polynomials (Hardy-Watson) \cite[Proposition 6.1.1]{Th}
asserts that: For $\a>-\frac12$ 
and $x, y\in \RR_+$,
\begin{align} \label{Hardy-Watson}
&\frac{\Gamma(n+1)}{\Gamma(n+\a+1)}L_n^\a(x^2)L_n^\a(y^2)\\
&\quad=\frac{2^\a}{\sqrt{2\pi}}\int_0^\pi L_n^\a
   \left(x^2+y^2+2xy\cos \theta\right)e^{-xy\cos \theta}
    j_{\a-1/2}(xy\sin \theta)\sin^{2\a}\theta d\theta,\notag
\end{align}
where $j_\a(x):=x^{-\a}J_\a(x)$ with $J_\a(x)$ being the Bessel function.

It will be convenient to denote $x^2 := (x_1^2, \ldots, x_d^2)$.
Combining (\ref{cesaro-kernel})-(\ref{Hardy-Watson}), we arrive at
\begin{align}\label{cesaroP}
P_n^{\a,\delta}(x^2,y^2)& = c_\a
\int_{[0,\pi]^d} P_n^{\a,\delta} (z(x,y,\theta),0) d\mu^\a_{x,y}(\theta)\\
& =  c_\a \int_{[0,\pi]^d} L_n^{|\a|+\delta+d}
\Big(\|x\|_2^2+\|y\|_2^2 + \sum_{i=1}^d x_iy_i \cos \theta_i\Big) d\mu^\a_{x,y}(\theta),\notag
\end{align}
where
$c_\a =  (2\pi)^{-d/2}2^{|\a|} \prod_{i=1}^d \Gamma(\a_i+1)$,
$z(x,y,\theta) = (z_1(x,y,\theta),\ldots,z_d(x,y,\theta))$ with
$z_i(x,y,\theta) = x_i^2+y_i^2 + 2 x_i y_i \cos \theta_i$, and
\begin{equation}\label{eq:dmu}
d \mu^{\a}_{x,y}(\theta) : = e^{-\sum_{i=1}^d x_iy_i\cos \theta_i}
\prod_{i=1}^d  j_{\a_i-1/2}(x_iy_i\sin \theta_i)\sin^{2\a_i}\theta_i d \theta.
\end{equation}

Some standard asymptotic properties of Laguerre functions will be needed.
The univariate Laguerre functions $\CL_n^\alpha$ are defined by
\begin{equation}\label{def.L-n}
\CL_n^\alpha(x):= \Big(\frac{\Gamma(n+1)}{\Gamma(n+\a+1)}\Big)^{1/2}
    e^{-x/2}x^{\a/2}L_n^\a(x).
\end{equation}


\begin{lem} \label{lem:asympt-Ln}
Set $N:=4n+2\alpha+2$. The Laguerre functions $\CL_n^\alpha$ satisfy
\begin{equation}\label{asympt-Ln}
|\CL_n^\a(x)| \le c\left\{
\begin{array}{ll}
(xN)^{\a/2}, & 0 < x\le 1/N,\\
(xN)^{-1/4}, & 1/N\le x \le N/2,\\
N^{-1/4}(N^{1/3}+|N-x|)^{-1/4}, &N/2 \le x \le 3N/2,\\
e^{-\gamma x}, & x\ge 3N/2,
\end{array}
\right.
\end{equation}
where $\gamma >0$ is an absolute constant.
\end{lem}

This lemma is contained in \cite[\S 8.22]{Sz} (see also
\cite[Lemma 1.5.3]{Th}).
Using that $\Gamma(n+\a+1)/\Gamma(n+1)
\sim n^\a$ one easily extracts from (\ref{asympt-Ln}) the estimates
\begin{equation}\label{est1-eLn}
e^{-x/2}|L_n^\a(x)| \le cn^{\alpha/2-1/4}x^{-\alpha/2-1/4},
\quad x\in \RR_+\setminus(N/2, 3N/2),
\end{equation}
and, for $N/2 \le x \le 3N/2$,
\begin{equation}\label{est4-eLn}
e^{-x/2}|L_n^\a(x)|
\le cx^{-\alpha/2}n^{\alpha/2-1/4}(n^{1/3}+|4n+2\a+2-x|)^{-1/4}.
\end{equation}
Also, from (\ref{asympt-Ln})
\begin{equation}\label{est2-eLn}
e^{-x/2}|L_n^\a(x)| \le cn^\alpha,
\quad x\in \RR_+,
\end{equation}
and since $\|\CL_n^\a\|_\infty \le c$, again by (\ref{asympt-Ln}),
\begin{equation}\label{est3-eLn}
e^{-x/2}|L_n^\a(x)| \le c(n/x)^{\a/2},
\quad x\in \RR_+.
\end{equation}


Let $K_n^\a(x,y)$ be the reproducing kernel of the Laguerre polynomials. Then
\begin{equation}\label{def-Kn}
K_n^\a(x,y)= c_\alpha\sum_{j=0}^n \frac{L_j^\a(x) L_j^\a(y)}{L_j^\a(0)},
\quad x,y\in \RR_+.
\end{equation}
The Christoffel function is defined by
\begin{equation} \label{def-Christoffel}
\lambda_n^\a(x) : = [K_n^\a(x,x)]^{-1}, \quad x \in \RR_+.
\end{equation}
For this function it is known that (see \cite{MO} and the references therein)
\begin{equation} \label{Christoph}
    c_1 \varphi_n(x) \le \frac{\lambda_n^\a(x)}{(x+\frac{1}{n})^\a e^{-x} }
          \le c_2 \varphi_n(x), \quad 0 \le x \le 4n,
\end{equation}
where
\begin{equation} \label{def-phi-n}
\varphi_n(x) := \sqrt{\frac{x + \frac1n}{4n-x+(4n)^{1/3}} }. 
\end{equation}

There are sharp estimates for $L_n^\a(x)$ in terms of $\varphi_n(x)$.
For any $x>0$, let $t_{k_x,n} $ denote the/a zero
of $L_n^\a(x)$ that is closest to $x$. Then (see e.g. \cite{MO})
\begin{equation}\label{lowerOP}
[L_n^\a(x)]^2 \left(x+\frac{1}{n}\right)^{\a+1} e^{-x}
\sim n^\a \varphi_n(x) \frac{(x-t_{k_x,n})^2}{(t_{k_x,n}- t_{k_x \pm 1,n})^2},
\quad x\in [t_{1,n}, t_{n,n}].
\end{equation}
Here and in the following
$t_{1,n}, \ldots, t_{n,n}$ denote the zeros of $L_n^\a(x)$.
They are known to satisfy \cite[\S 6.31]{Sz}
\begin{equation}\label{eq:zeros}
   c n^{-1} \le t_{1,n} < t_{2,n} < \ldots < t_{n,n} \le 4n + 2\a + 2 - c (4n)^{1/3}.
\end{equation}
Furthermore (see \cite[(6.31.11)]{Sz}),
\begin{equation}\label{zero-est}
c_*\frac{\nu^2}{n} \le t_{\nu,n} \le \frac{4\nu^2}{n}+ c(\a)\frac{\nu}{n}
\quad\mbox{and hence}\quad t_{\nu,n}\sim \frac{\nu^2}{n}.
\end{equation}
In addition (see \cite{MO} and the references therein),
\begin{equation}\label{zeros-difference}
t_{\nu+1,n} - t_{\nu,n} \sim \varphi_n(t_{\nu, n}).
\end{equation}
Therefore, if $\nu\le (1-\eps)n$ for some $\eps>0$, then by (\ref{zero-est})
$
t_{\nu,n} \le (1-\eps)^24n+c(\a),
$
and hence, using (\ref{zeros-difference}) and (\ref{def-phi-n}),
\begin{equation}\label{zeros-diff2}
t_{\nu+1,n} - t_{\nu,n} \sim \frac{\nu}{n}
\quad\mbox{if}\quad \nu\le (1-\eps)n.
\end{equation}
On the other hand, by (\ref{zeros-difference}) and (\ref{def-phi-n}), in general,
\begin{equation}\label{zeros-diff3}
\frac{c'}{n} \le t_{\nu+1,n} - t_{\nu,n} \le c''n^{1/3}.
\end{equation}

We will need the Gaussian quadrature formula with weight $t^\alpha e^{-t}$
on $(0, \infty)$ \cite{Sz}:
\begin{equation}\label{Gauss-0}
\int_0^\infty f(t)t^\a e^{-t}dt\sim \sum_{\nu=1}^n w_{\nu,n} f(t_{\nu,n}),
   \qquad w_{\nu,n}: = \lambda_n^\a(t_{\nu,n}),
\end{equation}
where $t_{\nu,n}$ are the zeros of $L_n^\a(t)$ and $\lambda_n^\a(x)$ is the
Christoffel function, defined in (\ref{def-Christoffel}).
This quadrature is exact for all algebraic polynomials of degree $2n-1$.

\section{Localized kernels associated with Laguerre functions}\label{localization}
\setcounter{equation}{0}

\subsection{The setting}

There are three kinds of univariate Laguerre functions considered in the literature
(see \cite{Th}), defined by
\begin{equation}\label{def.Psi-n}
\PsiP_n^\alpha(x):=
 \Big(\frac{2\Gamma(n+1)}{\Gamma(n+\a+1)}\Big)^{1/2}e^{-x^2/2}L_n^\a(x^2),
\end{equation}
$\CL_n^\alpha(x)$ have already been defined in (\ref{def.L-n}), and
\begin{equation}\label{def.M-n}
\CM_n^\alpha(x):= (2x)^{1/2}\CL_n^\alpha(x^2).
\end{equation}
It is well known that $\{\PsiP_n^\alpha\}_{n\ge 0}$ is an orthonormal basis
for the weighed space $L^2(\RR_+, x^{2\a+1})$,
while $\{\CL_n^\alpha\}_{n\ge 0}$ and $\{\CM_n^\alpha\}_{n\ge 0}$
are orthonormal bases for $L^2(\RR_+)$. 

Throughout this paper we will use standard multi-index notation.
Thus, for $x \in \RR^d$ and $\a \in \RR_+^d$, we write
$x^\a := x_1^{\a_1} \ldots x_d^{\a_d}$.
We will use ${\bf 1}$ to denote the vector ${\bf 1} := (1,1,\ldots,1)$.
Then, for instance, $x^{{\bf 1}/2} := x_1^{1/2} \cdots x_d^{1/2}$.
The $d$-dimensional Laguerre functions are defined by
\begin{align}
\PsiP_\nu^\a (x) &: = \PsiP_{\nu_1}^{\a_1}(x_1) \ldots \PsiP_{\nu_d}^{\a_d}(x_d),\label{def.dPsi_n}\\
\CL_\nu^\a (x) &: = \CL_{\nu_1}^{\a_1}(x_1) \ldots \CL_{\nu_d}^{\a_d}(x_d), \label{def.dL_n}\\
\CM_\nu^\a (x) &: = \CM_{\nu_1}^{\a_1}(x_1) \ldots \CM_{\nu_d}^{\a_d}(x_d), \label{def.dM_n}
\end{align}
where $\nu = (\nu_1,\ldots,\nu_d) \in \NN_0^d$ and $\a = (\a_1,\ldots, \a_d)$.
Clearly, $x^{-\a} e^{|x|} \CL_\nu^\a(x)$ is a polynomial of
degree $n=|\nu|=\nu_1+\dots+\nu_d$ and
\begin{equation} \label{rel-Psi-L}
\PsiP_\nu^\a (x) = 2^{d/2}x^{-\a} \CL_\nu^\a(x_1^2,\ldots,x_d^2).
\end{equation}
Evidently, $\{\PsiP_\nu^\alpha\}$ is an orthonormal basis
for the weighed space $L^2(\RR_+^d, \w)$, $\w(x):=x^{2\a+\1}$,
while $\{\CL_\nu^\alpha\}$ and $\{\CM_\nu^\alpha\}$
are orthonormal bases for $L^2(\RR_+^d)$ (with weight $1$).

We will utilize the basis $\{\cF^\a_\nu\}$ to the construction of frames for
the space $L^2(\w):=L^2(\RR_+^d, \w)$.
The same scheme based on $\{\CL_\nu^\alpha\}$ or $\{\CM_\nu^\alpha\}$
can be used for the construction of frames in $L^2(\RR_+^d)$.

As explained in the introduction, kernels of type (\ref{defn-Lambda-n})
will play a critical role in the present paper.
For our purposes we will be considering cut-off functions $\ha$
that satisfy:


\begin{defn} \label{defn:admissible}
A function $\ha \in C^\infty[0, \infty)$ is said to be admissible
of type $(a)$ or type $(b)$ if $\ha$ satisfies one of the following conditions:
\begin{enumerate}
\item[(a)]
$\supp \ha \subset [0,1+v]$, $\wh a(t) =1$ on $[0, 1]$, $v>0$; or
\item[(b)]
$\supp \ha \subset [u, 1+v]$, where $0<u<1$ and $v>0$.
\end{enumerate}
Here $u$, $v$ are fixed constants.
\end{defn}

\noindent
For an admissible function $\ha$ we introduce the kernels
\begin{align}
\LL_n(x,y)&:= \sum_{m=0}^\infty \ha\Big(\frac{m}{n}\Big) \PsiP_m^\a(x,y)
\quad \mbox{with}\quad
  \PsiP_m^\a(x,y):= \sum_{|\nu| =m}\PsiP_\nu^\a(x)\PsiP_\nu^\a(y),\label{def.Ln}\\
\tLL_n(x,y)&:= \sum_{m=0}^\infty \ha\Big(\frac{m}{n}\Big)  \CL_m^\a(x,y)
\quad \mbox{with}\quad
\CL_m^\a(x,y): = \sum_{|\nu|=m} \CL_\nu^\a(x)\CL_\nu^\a(y), \label{def.LLn}\\
\LL^*_n(x,y)&:= \sum_{m=0}^\infty \ha\Big(\frac{m}{n}\Big)  \CM_m^\a(x,y)
\quad \mbox{with}\quad
\CM_m^\a(x,y): = \sum_{|\nu|=m} \CM_\nu^\a(x)\CM_\nu^\a(y). \label{def.LLLn}
\end{align}
%
%
The rapid decay of the kernels
$\LL_n(x,y)$, $\tLL_n(x,y)$, and $\LL^*_n(x,y)$
and their partial derivatives
away from the main diagonal $y = x$ in $\RR_+^d \times \RR_+^d$
will be vital for our further development.

\subsection{The localization of \boldmath $\LL_n$ and its partial derivatives}

Recall the definition of the weight
$ 
\WW(n ;x):= \prod_{i=1}^d(x_i+n^{-\frac12})^{2\a_i+1}.
$ 


\begin{thm} \label{thm:est-Lambda-n}
Let $\wh a$ be admissible and let $\sigma > 0$. Then there is a constant $c_\sigma$
depending only on $\sigma$, $\alpha$, and $\ha$ such that for $x, y \in \RR_+^d$
\begin{equation} \label{est-Lambda-n}
|\LL_n(x,y)|
\le c_\sigma \frac{n^{d/2}}{\sqrt{\WW(n; x)}\sqrt{\WW(n; y)}(1+ n^{1/2}\|x-y\|)^\sigma},
\end{equation}
and furthermore, for $1 \le r \le d$,
\begin{equation}\label{est-derivative}
 \left| \frac{\partial}{\partial x_r } \LL_n(x,y) \right|
 \le c_\sigma\frac{ n^{(d+1)/2}}
  {\sqrt{\WW(n; x)}\sqrt{\WW(n; y)}(1+ n^{1/2}\|x-y\|)^\sigma} .
\end{equation}
Here the dependence of $c_\sigma$ on $\ha$ is of the form
$c_\sigma=c(\sigma, \alpha)\max_{0\le l \le k}\|\ha^{(l)}\|_{L^\infty}$,
where $k \ge \sigma +2|\alpha|+d/2$.

In addition to this, there exists a constant $\expon>0$ such that
if $x, y \in \RR_+^d$ and
$\max\{\|x\|, \|y\|\} \ge (6(1+v)n+3\alpha+3)^{1/2}$, then
\begin{equation} \label{est-Lambda-n-large}
|\LL_n(x,y)|
\le c_\sigma \frac{e^{-\expon \max\{\|x\|, \|y\|\}^2}}{(1+ n^{1/2}\|x-y\|)^\sigma}
\end{equation}
and, for $1 \le r \le d$,
\begin{equation}\label{est-derivative-large}
 \left| \frac{\partial}{\partial x_r } \LL_n(x,y) \right|
 \le c_\sigma\frac{e^{-\expon \max\{\|x\|, \|y\|\}^2}}
  {(1+ n^{1/2}\|x-y\|)^\sigma}.
\end{equation}

\end{thm}

To keep our exposition more fluid we relegate the proofs of these
and the estimates to follow in this section to \S~\ref{Proofs}.

We next use estimate (\ref{est-Lambda-n}) to bound the $L^p$-integral of
$\LL_n(x, y)$, in particular, we show that
$\int_{\RR_+^d} |\LL_n(x, y)|\w(y) dy \le c <\infty$.


\begin{prop} \label{prop:Lp-bound}
For $0 < p < \infty$, we have
\begin{equation}\label{Lp-bound}
\int_{\RR^d_+} |\LL_n(x,y)|^p \w(y)dy
\le cn^{(d/2)(p-1)}\WW(n; x)^{-(p-1)},
\quad x \in \RR_+^d.
\end{equation}
\end{prop}

Estimate (\ref{Lp-bound}) is immediate from (\ref{est-Lambda-n}) and the following lemma
which will be instrumental in the subsequent development.


\begin{lem}\label{lem:instrument}
If $s\in \RR$ and $\sigma > d((2\|\alpha\|+1)(|s|+1)+1)$, then
\begin{equation} \label{instrument}
\int_{\RR^d_+} \frac{\w(y)\, dy}{\WW(n; y)^s(1+n^{1/2}\|x-y\|)^\sigma}
\le \frac{cn^{-d/2}}{\WW(n;x)^{s-1}},
\quad x\in\RR^d_+.
\end{equation}
\end{lem}


We next give a lower bound estimate:


\begin{thm}\label{thm:lower-bound}
Let $\ha$ be admissible in the sense of Definition~\ref{defn:admissible}
and $|\ha|>c_\diamond>0$ on $[1, 1+\tau]$, $\tau>0$.
Then for any $\delta>0$
\begin{equation}\label{lower-bound}
\int_{\RR_+^d} |\LL_n(x,y)|^2 \w(y)dy
\ge c \,n^{d/2}\WW(n; x)^{-1}
\quad x \in [0,\sqrt{(4-\delta)n}]^d,
\end{equation}
where $c>0$ depends only on $\a$, $d$, $\tau$, $\delta$, and $c_\diamond$.
\end{thm}

By the orthogonality of the Laguerre functions it readily follows that
$$
\int_{\RR^d_+}    |\LL_n(x,y)|^2 \w(y)dy
= \sum_{m=0}^\infty |\ha(m/n)|^2 \PsiP_m^\a(x,x),
$$
and hence Theorem~\ref{thm:lower-bound} is an immediate consequence of the
following lemma.


\begin{lem}\label{lem:lower-bound}
For any $\eps>0$ and $\delta>0$ there exists a constant $c > 0$ such that
\begin{equation} \label{lowerbd}
\sum_{m=n}^{n + \lfloor d\varepsilon n\rfloor}
\PsiP_m^\a (x,x) \ge c n^{d/2}\WW(n; x)^{-1},
\quad x \in [0,\sqrt{(4-\delta)n}]^d.
\end{equation}
\end{lem}

\subsection{The localization of \boldmath $\tLL_n$ and its partial derivatives}
The localization of the kernels $\tLL_n$ can be deduced from
the localization of $\LL_n$ given above.


\begin{thm} \label{thm:est-tilde-Lambda-n}
Let $\wh a$ be admissible. Then for any $\sigma > 0$ there is a constant $c_\sigma>0$
depending only on $\sigma$, $\alpha$, and $\ha$ such that for $ x, y \in \RR_+^d$,
\begin{equation} \label{est-tilde-Lambda-n}
|\tLL_n(x,y)|
\le c_\sigma \frac{n^{d/2}}{\prod_{i=1}^d(x_i+n^{-1})^{\frac14}(y_i+n^{-1})^{\frac14}
      (1+ n^{1/2} \|x^{\1/2}- y^{\1/2}\|)^\sigma},
\end{equation}
and, for $1 \le r \le d$,
\begin{align}\label{est-tilde-derivative}
& \left| \frac{\partial}{\partial x_r }\wt \LL_n(x,y) \right|
 \le \frac{c n^{d/2+1}}
  {\prod_{i=1}^d(x_i+n^{-1})^{\frac14}(y_i+n^{-1})^{\frac14}
      (1+ n^{1/2}\|x^{\1/2}- y^{\1/2}\|)^\sigma}.
\end{align}
Here the dependence of $c_\sigma$ on $\ha$ is as in Theorem~\ref{thm:est-tilde-Lambda-n}.
\end{thm}

Estimates for $\tLL_n$ like the ones of (\ref{est-Lambda-n-large})-(\ref{lower-bound})
can be extracted from (\ref{est-Lambda-n-large})-(\ref{lower-bound}).
The results from this and the next subsections follow easily
from Theorem~\ref{thm:est-Lambda-n}, see \S\ref{Other-loc-est}.

\subsection{The localization of \boldmath $\LL_n^*$ and its partial derivatives}
The localization properties of $\LL_n^*(x,y)$ appear simpler:


\begin{thm} \label{thm:est-Lambda*-n}
Let $\wh a$ be admissible. Then for any $\sigma > 0$ there is a constant $c_\sigma$
such that for $x, y \in \RR_+^d$
\begin{equation} \label{est-Lambda*-n}
|\LL_n^*(x,y)|
\le c_\sigma \frac{n^{d/2}}{(1+ n^{1/2}\|x-y\|)^\sigma},
\end{equation}
and, for $1 \le r \le d$,
\begin{equation}\label{est-derivative*}
 \left| \frac{\partial}{\partial x_r } \LL_n^*(x,y) \right|
 \le c_\sigma\frac{ n^{(d+1)/2}}
  {(1+ n^{1/2}\|x-y\|)^\sigma}.
\end{equation}
\end{thm}

Estimates for $\LL_n^*$ similar to the ones of (\ref{est-Lambda-n-large})-(\ref{lower-bound})
can easily be obtained.


\section{Additional background material}\label{background}
\setcounter{equation}{0}


\subsection{Norm equivalence}\label{norms}


\begin{prop}\label{prop:norms}
Let $0<q\le p\le \infty$ and $g\in V_n$ $(n\ge 1)$. Then
\begin{equation} \label{norms1}
\|g\|_p\le cn^{(d+|\alpha|)(1/q-1/p)}\|g\|_q
\end{equation}
and, for any $s\in\RR$,
\begin{equation} \label{norms2}
\|\WW(n;\cdot)^s g(\cdot)\|_p
\le cn^{(d/2)(1/q-1/p)}\|\WW(n;\cdot)^{s+1/p-1/q} g(\cdot)\|_q.
\end{equation}
Furthermore, for any $s\in\RR$
\begin{equation} \label{norms3}
\|g\|_p
\le cn^M\|\WW(n;\cdot)^s g(\cdot)\|_q,
\end{equation}
where
$M$ depends only on $\a, d, p, q$, and $s$.
\end{prop}
The proof of this proposition employs the localized kernels from \S\ref{localization}
and is rather standard.
For completeness we give it in \S\ref{Proofs}.

\subsection{Maximal operator}\label{max-operator}

We define the ``cube" centered at $\xi\in\RR^d_+$ of ``radius" $r>0$ by
$Q_\xi(r):=\{x\in\RR^d_+: \|x-\xi\|<r\}$.
Let $\cM_t$ be the maximal operator, defined by
\begin{equation}\label{maxOperator}
\cM_t f(x) := \sup_{Q:\, x \in Q}
\left(\frac{1}{\mu(Q)} \int_Q |f(y)|^t\w(y)\,dy \right)^{1/t},
\quad x\in \RR^d_+,
\end{equation}
where the $\sup$ is over all ``cubes" $Q$ in $\RR^d_+$ with sides parallel
to the coordinate axes which contain $x$.
It is easy to see that
\begin{equation}\label{size-muQ}
\mu(Q_\xi(r))\sim r^d\prod_{j=1}^d (\xi_j+r)^{2\alpha_j+1}.
\end{equation}
Hence
$\mu(Q_\xi(2r)) \le c\mu(Q_\xi(r))$,
i.e. $\mu(\cdot)$ is a doubling measure.
Therefore, the theory of maximal operators applies and
the Fefferman-Stein vector-valued maximal inequality is valid
(see \cite{Stein}):
If $0<p<\infty$, $0<q\le \infty$, and $0 < t < \min \{p, q\}$,
then for any sequence of functions  $f_1, f_2, \dots$ on $\RR^d_+$
\begin{equation}\label{max-ineq}
\Big\|\Big( \sum_{j=1}^{\infty}
     \left[\cM_t f_j(\cdot)\right]^q \Big)^{1/q}\Big\|_p
\le c \Big\|\Big( \sum_{j=1}^{\infty}|f_j(\cdot)|^q \Big)^{1/q}\Big\|_p,
\end{equation}
where $c = c(p,q,t,d,\alpha)$.

\subsection{Distributions on \boldmath $\RR^d_+$}

We will use as test functions the set $\cS_+$ of all functions
$\phi\in C^\infty([0, \infty)^d)$
such that
\begin{equation}\label{schwartz}
P_{\b,\g}(\phi):=\sup_{x\in \RR^d_+}|x^\g \DD^\b \phi(x)| <\infty
\quad\mbox{for all multi-idices $\g$ and $\b$,}
\end{equation}
with the topology on $\cS_+$ defined by the semi-norms $P_{\b,\g}$.
Then the space $\cS_+'$ of all temperate distributions on $\RR^d_+$
is defined as the set of all continuous linear functionals on $\cS_+$.
The pairing of $f\in \cS_+'$ and $\phi\in\cS_+$ will be denoted by
$\langle f, \phi \rangle := f(\overline{\phi})$ which is consistent
with the inner product
$\langle f, g \rangle := \int_{\RR^d_+}f(x) \overline{g(x)}\w(x)dx$
in $\L2(\RR^d_+, \w)$.

It will be convenient for us to introduce the following
``convolution":


\begin{defn}\label{def:convolution}
For functions $\Phi: \RR^d_+\times\RR^d_+ \to \CC$ and $f: \RR^d_+ \to \CC$,
we define
\begin{equation}\label{convolution}
\Phi*f(x) := \int_{\RR^d_+} \Phi(x, y)f(y)\w(y)dy.
\end{equation}
In general,
if $f \in \cS_+'$ and $\Phi: \RR^d_+\times\RR^d_+\to\CC$ is such that
$\Phi(x, y)$ belongs to $\cS_+$ as a function of $y$ $(\Phi(x, \cdot)\in \cS_+)$,
we define $\Phi*f$ by
\begin{equation}\label{convolution1}
\Phi*f (x) := \langle f, \overline{\Phi(x, \cdot)} \rangle,
\end{equation}
where on the right $f$ acts on $\overline{\Phi(x, y)}$ as a function of $y$.
\end{defn}

We now give some properties of the above convolution that can be proved
in a standard way.


\begin{lem}\label{lem:convolution}
$(a)$ If $f \in \cS_+'$ and $\Phi (\cdot, \cdot) \in \cS_+(\RR^d_+\times \RR^d_+)$, then
$\Phi*f \in \cS_+$.
Furthermore $\cF^\a_n*f\in \V_n$.

$(b)$ If $f \in \cS_+'$, $\Phi (\cdot, \cdot) \in \cS_+(\RR^d_+\times \RR^d_+)$, and
$\phi \in \cS_+$, then
$\langle \Phi*f, \phi \rangle = \langle f, \overline{\Phi}*\phi \rangle$.

$(c)$ If $f \in \cS_+'$,
$\Phi (\cdot, \cdot), \Psi (\cdot, \cdot) \in \cS_+(\RR^d_+\times \RR^d_+)$,
and
$\Phi (y, x)= \Phi (x, y)$, $\Psi (y, x)= \Psi (x, y)$,
then
\begin{equation}\label{Psi*Phi*f}
\Psi*\overline{\Phi}*f(x)
= \langle \Psi(x, \cdot), \Phi(\cdot, \cdot)\rangle*f.
\end{equation}
\end{lem}

Evidently the Laguerre functions $\{\cF^\a_\nu\}$ belong to $\cS_+$.
Moreover, the functions in $\cS_+$ can be
characterized by the coefficients in their Laguerre expansions.
Denote
\begin{equation}\label{def-P*}
P_r^{*}(\phi) := \sum_{n=0}^\infty (n+1)^r\|\cF^\a_n*\phi\|_{2}
=\sum_{n=0}^\infty (n+1)^r
\Big(\sum_{|\nu|=n}|\langle \phi, \cF^\a_\nu\rangle|^2\Big)^{1/2}.
\end{equation}


\begin{lem}\label{lem:char-S} A function
$\phi\in\cS_+$ if and only if
$|\langle \phi, \cF^\a_\nu\rangle| \le c_k(|\nu|+1)^{-k}$
for all multi-indices $\nu$ and all $k$.
Moreover, the topology in $\cS_+$ can be equivalently defined by
the semi-norms $P_r^{*}$.
\end{lem}
The proof of this lemma is given in \S\ref{Proofs}.


\section{Construction of frame elements (Needlets)}\label{def-needlets}
\setcounter{equation}{0}

In this section we construct frames utilizing the localized kernels
from \S\ref{localization} and a~cubature formula on $\RR_+^d$.
As explained in the introduction, we will only use the Laguerre functions
$\{\PsiP_\nu^\a \}$ defined in (\ref{def.dPsi_n}).

\subsection{Cubature formula}\label{cubature}

We will utilize the Gaussian quadrature (\ref{Gauss-0}) for the construction
of the needed cubature formula on $\RR_+^d$.
Given $n\ge 1$, we define, for $\nu=1, \dots, n$,
\begin{equation}\label{def-qub1}
\xi_{\nu, n}:=\sqrt{t_{\nu, n}}
\quad\mbox{and}\quad
c_{\nu, n}:= \frac12 w_{\nu, n}e^{t_{\nu, n}}
=\frac12 \lambda_n^\alpha(t_{\nu, n})e^{t_{\nu, n}}
=\frac12 \lambda_n^\alpha(\xi_{\nu, n}^2)e^{\xi_{\nu, n}^2},
\end{equation}
where $\{t_{\nu, n}\}$ are the zeros of $L_n^\a(t)$ and
$\{w_{\nu, n}\}$ are the weights from (\ref{Gauss-0}).

It follows by (\ref{zero-est}) and (\ref{zeros-diff2})-(\ref{zeros-diff3})
that
\begin{equation}\label{zero4}
\xi_{\nu, n}\sim \frac{\nu}{\sqrt{n}},
\end{equation}
\begin{equation}\label{zero-diff4}
\xi_{\nu+1, n}- \xi_{\nu, n} \sim n^{-1/2}
\quad\mbox{if}\quad
1\le \nu\le (1-\eps)n,
\end{equation}
and, in general,
\begin{equation}\label{zero-diff5}
c_1n^{-1/2} \le \xi_{\nu+1, n}- \xi_{\nu, n} \le c_2n^{-1/6}.
\end{equation}
Furthermore, using (\ref{Christoph}) and (\ref{zeros-difference})
we obtain
\begin{equation}\label{est-lambda}
c_{\nu, n}\sim \ph_n(t_{\nu, n})t_{\nu, n}^\a
\sim (t_{\nu+1, n}-t_{\nu, n})t_{\nu, n}^\a
\sim (\xi_{\nu+1, n}-\xi_{\nu, n})\xi_{\nu, n}^{2\a+1}.
\end{equation}

Now, for $\gamma=(\gamma_1, \dots, \gamma_d)\in\NN_0^d$
we set
\begin{equation}\label{def-qub2}
c_{\gamma, n}:=\prod_{j=1}^d c_{\gamma_j, n}
\quad\mbox{and}\quad
\xi_{\gamma, n}:=(\xi_{\gamma_1, n}, \dots, \xi_{\gamma_d, n}).
\end{equation}


\begin{prop}\label{prop:cubature}
The cubature formula
\begin{equation}\label{cubature1}
\int_{\RR_+^d}f(x)g(x)\w(x)dx
\sim \sum_{\gamma_1=1}^n \cdots\sum_{\gamma_d=1}^n
                c_{\gamma, n}f(\xi_{\gamma, n})g(\xi_{\gamma, n})
\end{equation}
is exact for all $f\in V_\ell$ and $g\in V_m$ provided $\ell+m\le 2n-1$.
\end{prop}

\noindent
{\bf Proof.}
Evidently, it suffices to consider only the case $d=1$.
Suppose $f\in V_\ell$ and $g\in V_m$ with $\ell+m\le 2n-1$.
Let
$f(x)=:F(x^2)e^{-x^2/2}$ and
$g(x)=:G(x^2)e^{-x^2/2}$,
where $F\in \Pi_\ell^1$, $G\in\Pi_m^1$
with $\Pi_j^1$ being the set of all univariate polynomials of degree $\le j$.
Then using the properties of quadrature formula (\ref{Gauss-0}), we get
\begin{align*}
\int_0^\infty f(x)g(x)\w(x)dx
&=\int_0^\infty F(x^2)G(x^2)x^{2\a+1}e^{-x^2}dx
=\frac12 \int_0^\infty F(t)G(t)t^{\a}e^{-t}dt\\
&=\frac12 \sum_{\nu=1}^n w_{\nu, n} F(t_{\nu, n})G(t_{\nu, n})
=\sum_{\nu=1}^n \frac12 w_{\nu, n} F(\xi_{\nu, n}^2)G(\xi_{\nu, n}^2)\\
&
=\sum_{\nu=1}^n \frac12 \lambda_n^\a(\xi_{\nu, n}^2)e^{\xi_{\nu, n}^2}
                                             f(\xi_{\nu, n})g(\xi_{\nu, n}),
\end{align*}
which completes the proof.
$\qed$

To construct our frame elements we need the cubature formula from (\ref{cubature1}) with
\begin{equation}\label{def-nj}
n=n_j:=\lfloor c_*^{-1}(1+11\delta)\sqrt{6}\cdot 4^j\rfloor +1\sim 4^j,
\end{equation}
where $0< c_* \le 1$ is the constant from (\ref{zero-est})
and $0<\delta<1/26$ is an arbitrary but fixed constant.
For $j\ge 0$, we define
\begin{equation}\label{def-Xj}
\cX_j:=\{\xi\in\RR_+^d: \xi
=\xi_{\gamma, n_j}, 1\le \gamma_\ell\le n_j, 1\le \ell\le d\}.
\end{equation}
Note that $\#\cX_j=n_j^d\sim 4^{jd}$.
Now, if $\xi\in\cX_j$ and $\xi=\xi_{\gamma, n_j}$,
we set $c_\xi:= c_{\gamma, n_j}$.

As an immediate consequence of Proposition~\ref{prop:cubature} we get


\begin{cor}\label{cor:cubature}
The cubature formula
\begin{equation}\label{cub-cor}
\int_{\RR_+^d}f(x)g(x)\w(x)dx
\sim \sum_{\xi\in\cX_j} c_\xi f(\xi)g(\xi)
\end{equation}
is exact for all $f\in V_\ell$ and $g\in V_m$ provided $\ell+m\le 2n_j-1$.
\end{cor}

\bigskip


\noindent
{\bf Tiling.} We next introduce rectangular tiles $\{R_\xi\}$ with ``centers"
at the points $\xi\in\cX_j$.
Set $I_1:=[0, (\xi_1+\xi_2)/2]$ and
$$
I_\nu:= [(\xi_{\nu-1}+\xi_\nu)/2, (\xi_{\nu}+\xi_{\nu+1})/2],
\quad \nu=2, \dots, n_j,
$$
where
$\xi_\nu:= \xi_{\nu, n_j}$ , $\nu=1, \dots, n_j$,
are from (\ref{def-qub1}) and $\xi_{n_j+1}:= \xi_{n_j}+2^{j/3}$.

To every $\xi=\xi_\gamma=(\xi_{\gamma_1}, \dots, \xi_{\gamma_d})$ in $\cX_j$
we associate a tile $R_\xi$ defined by
\begin{equation}\label{def-Rxi}
R_\xi:= I_{\gamma_1}\times \cdots \times I_{\gamma_d}.
\end{equation}
We also set
\begin{equation}\label{def-Qj}
Q_j:=\cup_{\xi\in\cX_j} R_\xi.
\end{equation}
Evidently, different tiles $R_\xi$ do not overlap and
$Q_j\sim [0, 2^j]^d$.

By (\ref{est-lambda}) it readily follows that
\begin{equation}\label{est-cn}
c_\xi \sim \mu(R_\xi):= \int_{R_\xi}\w(x)dx
\sim |R_\xi|\w(\xi) \sim |R_\xi|\WW(4^j;\xi).
\end{equation}
Assume $\xi\in\cX_j$, $\xi:=\xi_\gamma$, and $\|\xi\|\le (1+4\delta)\sqrt{6}\cdot 2^{j}$.
By (\ref{zero-est}) $\|\xi_\gamma\| \ge c_*^{1/2}\|\gamma\|n_j^{-1/2}$ and hence
$\|\gamma\|\le c_*^{-1/2}(1+4\delta)\sqrt{6}\cdot 2^{j}n_j^{1/2}\le (1-\delta)n_j$,
where the last inequality follows by the selection of $n_j$ in (\ref{def-nj}).
Therefore, for $\xi\in\cX_j$
\begin{equation}\label{size-R-xi1}
R_\xi\sim \xi+[-2^{-j}, 2^{-j}]^d
\;\;\hbox{and}\;\; \mu(R_\xi)\sim 2^{-jd}\w(\xi)
\;\;\hbox{if \; $\|\xi\|\le (1+4\delta)\sqrt{6}\cdot 2^{j}$,}
\end{equation}
while in general, for some positive constants $c_1, c_2, c', c''$,
\begin{equation}\label{size-R-xi2}
\xi+[-c_12^{-j}, c_12^{-j}]^d \subset R_\xi\subset \xi+[-c_22^{-j/3}, c_22^{-j/3}]^d
\quad \mbox{and}
\end{equation}
\begin{equation}\label{size-R-xi3}
c'2^{-jd}\w(\xi) \le \mu(R_\xi)\le c''2^{-jd/3}\w(\xi).
\end{equation}

The following simple inequality is immediate from the definition of
$\WW(n; x)$ in (\ref{def-Wn}) and will be useful in what follows:
\begin{equation}\label{W<W}
\WW(4^j; y)\le \WW(4^j; x)(1+2^j\|x-y\|)^{2|\a|+d},
\quad x, y\in \RR^d_+.
\end{equation}

\subsection{Definition of Needlets}\label{Def-of-needlets}

Let $\ha$, $\hb$
satisfy the conditions:  
\begin{equation}\label{ab1}
\ha, \hb \in C^\infty(\R),
\quad\supp \ha, \supp \hb \subset [1/4, 4],
\end{equation}
\begin{equation}\label{ab3}
|\ha(t)|, |\hb(t)|>c>0 \quad \mbox{if}~~
t \in [1/3, 3],
\end{equation}
\begin{equation}\label{ab5}
\overline{\ha(t)}\;\hb(t) + \overline{\ha(4t)}\;\hb(4t) =1
\quad\mbox{if}~~ t \in [1/4, 1].
\end{equation}
Hence,
\begin{equation}\label{unity1}
\sum_{m=0}^\infty \overline{\ha(4^{-m}t)}\;\hb(4^{-m}t) = 1,
\quad t \in [1, \infty).
\end{equation}

It is readily seen that (e.g. \cite{F-J2})
for any $\ha$ satisfying (\ref{ab1})-(\ref{ab3})
there exists $\hb$ satisfying (\ref{ab1})-(\ref{ab3})
such that $(\ref{ab5})$ holds.

\smallskip

Let $\ha$, $\hb$ satisfy (\ref{ab1})-(\ref{ab5}).
Then we set
\begin{align}
&\Phi_0(x, y) :=\PsiP_0^\a(x, y),
\quad
\Phi_j(x, y) := \sum_{m=0}^\infty
\ha\Big(\frac{m}{4^{j-1}}\Big)\PsiP_m^\a(x, y),
\quad \mbox{and}  \label{def.Phi-j}\\
&\Psi_0(x, y) :=\PsiP_0^\a(x, y),
\quad
\Psi_j(x, y) := \sum_{m=0}^\infty
\hb\Big(\frac{m}{4^{j-1}}\Big)\PsiP_m^\a(x, y),
\quad j\ge 1. \label{def-Psi-j}
\end{align}
Let $\cX_j$ be the set defined in (\ref{def-Xj})
and let $c_\xi$ be the coefficients of cubature formula (\ref{cub-cor}).
We define the $j$th level {\em needlets} by
\begin{equation}\label{def-needlets1}
\ph_\xi(x) := c_\xi^{1/2}\Phi_j(x, \xi)
\quad\mbox{and}\quad
\psi_\xi(x) := c_\xi^{1/2}\Psi_j(x, \xi),
\qquad \xi \in \cX_j.
\end{equation}
Set $\cX := \cup_{j = 0}^\infty \cX_j$.
We will use $\cX$ as an index set for our needlet systems
$\Phi$ and $\Psi$.
For this reason, (possibly) identical points from different levels $\cX_j$ are considered
as distinct elements of $\cX$.
We define
\begin{equation}\label{def-needlets2}
\Phi:=\{\ph_\xi\}_{\xi\in\cX}, \quad \Psi:=\{\psi_\xi\}_{\xi\in\cX}.
\end{equation}
We will term
$\{\ph_\xi\}$ {\em analysis needlets}
and $\{\psi_\xi\}$ {\em synthesis needlets}.

\medskip

\noindent
{\bf Localization of Needlets.}
An immediate consequence of Theorem~\ref{thm:est-Lambda-n} is the estimate:
For any $\sigma >0$ there exists a constant $c_\sigma >0$ such that
for all $x, y\in\RR^d_+$
\begin{equation}\label{local-Needlets1}
|\Phi_j(x, y)|, |\Psi_j(x, y)|
\le \frac{c_\sigma 2^{jd}}{\sqrt{\WW(4^j, x)}\sqrt{\WW(4^j, y)}(1+2^j\|x-y\|)^\sigma},
\end{equation}
while $c_\sigma 2^{jd}$ can be replaced by $c(\sigma, L) 2^{-jL}$
if $\max\{\|x\|, \|y\|\}\ge (1+\delta)\sqrt{6}\cdot 2^{j}$,
where $L>0$ is an arbitrary constant but the constant $c(\sigma, L)$
depends on $L$ as well.
We employ (\ref{local-Needlets1}) and (\ref{est-cn})
to obtain for $\xi\in\cX_j$
\begin{equation}\label{local-needlets2}
|\ph_\xi(x)|, |\psi_\xi(x)|
\le \frac{c 2^{jd/2}}{\sqrt{\WW(4^j, x)}(1+2^j\|x-\xi\|)^\sigma},
\;\;\; x\in\RR^d_+,
\end{equation}
and
\begin{equation}\label{local-needlets3}
|\ph_\xi(x)|, |\psi_\xi(x)|
\le \frac{c2^{-jL}}{\sqrt{\WW(4^j, x)}(1+2^j\|x-\xi\|)^\sigma},
\; \mbox{if \; $\|\xi\|\ge (1+\delta)\sqrt{6}\cdot 2^{j}$.}
\end{equation}

We next show that $\cS'_+$ and $L^p(\RR^d_+)$ have discrete
decompositions via needlets.


\begin{prop}\label{prop:needlet-rep}
$(a)$ If $f \in \cS'_+$, then
\begin{align}
f &= \sum_{j=0}^\infty
\Psi_j*\overline{\Phi}_j*f
\quad\mbox{in} \;\; \cS'_+ \;\; \mbox{and} \label{Needle-rep}\\
f &= \sum_{\xi \in \cX}
\langle f, \ph_\xi\rangle \psi_\xi
\quad\mbox{in} \;\; \cS'_+.\label{needlet-rep}
\end{align}

$(b)$ If $f \in L^p(\w)$, $1\le p < \infty$, then
$(\ref{Needle-rep})-(\ref{needlet-rep})$ hold in $L^p(\w)$.
Moreover, if $1 < p < \infty$, then the convergence in
$(\ref{Needle-rep})-(\ref{needlet-rep})$ is unconditional.
\end{prop}

\noindent
{\bf Proof.}
(a)
Note that $\Psi_j*\overline{\Phi}_j(x, y)$ is well defined
since $\Psi_j(x, y)$ and $\Phi_j(x, y)$
are symmetric functions (e.g. $\Psi_j(y, x)=\Psi_j(x, y)$).
By (\ref{def.Phi-j})-(\ref{def-Psi-j}) it follows that
$\Psi_0*\overline{\Phi}_0=\PP_0$
and
\begin{equation}\label{rep-Psi*Phi-1}
\Psi_j*\overline{\Phi}_j(x, y)
=\sum_{m=4^{j-2}}^{4^j} \overline{\ha\Big(\frac{m}{4^{j-1}}\Big)}
\hb\Big(\frac{m}{4^{j-1}}\Big)\PP_m(x, y),
\quad j\ge 1.
\end{equation}
Hence, (\ref{unity1}) and Lemma~\ref{lem:char-S} imply (\ref{Needle-rep}).
Evidently,
$\Psi_j(x, \cdot)$ and $\overline{\Phi_j(y, \cdot)}$
belong to $\V_{4^j}$ and
using the cubature formula from Corollary~\ref{cor:cubature},
we infer
\begin{align*}
\Psi_j*\overline{\Phi}_j(x, y)
&= \int_{\RR^d_+} \Psi_j(x, u)\overline{\Phi_j(y, u)}\,du\\
&=\sum_{\xi\in \cX_j}
c_\xi\Psi_j(x, \xi)\overline{\Phi_j(y, \xi)}
= \sum_{\xi\in \cX_j}\psi_\xi(x)\overline{\ph_\xi(y)}.\notag
\end{align*}
Therefore,
$
\Psi_j*\overline{\Phi}_j*f = \sum_{\xi\in \cX_j}
\langle f,\ph_\xi\rangle \psi_\xi
$
and combining this with (\ref{Needle-rep}) gives (\ref{needlet-rep}).

(b) In $L^p$ identity (\ref{Needle-rep}) follows easily by the rapid decay
of the kernels of the $n$th partial sums. We skip the details.
In $L^p$, identity (\ref{needlet-rep}) follows as above.
The unconditional convergence in $L^p(\w)$, $1<p<\infty$, is a consequence of
Proposition~\ref{prop:identification} and Theorem~\ref{thm:needlet-Tr-Liz} below.
$\qed$


\begin{re}\label{rem:frame}
Suppose that in the needlet construction
$\hb = \ha$ and $\ha\ge 0$.
Then
$\ph_\xi = \psi_\xi$
and $(\ref{needlet-rep})$ becomes
$f = \sum_{\xi \in \cX} \langle f, \psi_\xi\rangle \psi_\xi$.
It is easily seen that this representation holds
in $\L2$ and
$\|f\|_{\L2} =
\Big(\sum_{\xi \in \cX}
|\langle f, \psi_\xi\rangle|^2\Big)^{1/2}$,
$f\in \L2,$
i.e.
$\{\psi_\xi\}_{\xi \in \cX}$ is a tight frame for $\L2(\RR^d_+, \w)$.
\end{re}

\section{Laguerre-Triebel-Lizorkin  spaces}\label{Tri-Liz}
\setcounter{equation}{0}

We follow the general idea of using spectral decompositions (see e.g. \cite{Peetre}, \cite{T1})
to introduce Triebel-Lizorkin spaces on $\RR^d_+$ in the context of Laruerre expansions.
Our main goal is to show that these spaces can be characterized via needlet representations.

\subsection{Definition of Laguerre-Triebel-Lizorkin  spaces}\label{def-F-spaces}

Let a sequence of kernels $\{\Phi_j\}$ be defined by
\begin{equation}\label{def-Phi-j}
\Phi_0(x, y) := \cF_0^\a(x, y)
\quad\mbox{and}\quad
\Phi_j(x, y) := \sum_{m=0}^\infty \ha
\Big(\frac{m}{4^{j-1}}\Big)\cF_m^\a(x,y), ~~~j\ge 1,
\end{equation}
where $\{\cF_m^\a(x,y)\}$ are from \eqref{def.Ln} and
$\ha$ obeys the conditions
\begin{align}
&\quad \ha\in C^\infty[0, \infty),
\quad \supp \, \ha \subset [1/4, 4], \label{ha1}\\
&\quad  |\ha(t)|>c>0, \quad \text{if } t \in [1/3, 3].\label{ha2}
\end{align}


\begin{defn}\label{defn-Tr-Liz}
Let  $s, \r \in \RR$, $0<p<\infty$, and $0<q\le\infty$. Then the
Laguerre-Triebel-Lizorkin space $\Fsrpq:=\Fsrpq(\cF^\a)$ is defined as
the set of all distributions $f\in \cS_+'$ such that
\begin{equation}\label{Tri-Liz-norm}
\|f\|_{\Fsrpq}:=\Big\|\Big(\sum_{j=0}^{\infty}
\Big[2^{sj}\WW(4^j;\cdot)^{-\r/d}|\Phi_j*f(\cdot)|\Big]^q\Big)^{1/q}\Big\|_p <\infty
\end{equation}
with the usual modification when $q=\infty$.
\end{defn}

As is shown in Theorem~\ref{thm:needlet-Tr-Liz} below
the above definition is independent of the choice
of $\ha$ as long as $\ha$ satisfies (\ref{ha1})-(\ref{ha2}).


\begin{prop}\label{Fsrpq-embedding}
For all $s,\r \in \RR$, $0<p<\infty$, and $0<q\le\infty$,
$\Fsrpq$ is a $($quasi-$)$Banach space which is continuously embedded in $\cS_+'$.
\end{prop}

\noindent
{\bf Proof.}
The completeness of the space $\Fsrpq$
follows easily (see e.g. \cite{T1}, p. 49) by the continuous
embedding of $\Fsrpq$ in $\cS_+'$, which we establish next.

Let $\{\Phi_j\}$ be the kernels from the definition of $\Fsrpq$ with $\ha$ obeying
(\ref{ha1})-(\ref{ha2}) that are the same as (\ref{ab1})-(\ref{ab3}).
As already indicated there exists a function $\hb$ satisfying (\ref{ab1})-(\ref{ab5}).
We use this function to define $\{\Psi_j\}$ as in (\ref{def-Psi-j}).
Assume $f\in\Fsrpq$.
Then by Proposition~\ref{prop:needlet-rep}
$f=\sum_{j=0}^\infty \Psi_j*\Phi_j*f$ in $\cS'_+$
and hence
$$
\langle f, \phi\rangle=\sum_{j=0}^\infty \langle \Psi_j*\overline{\Phi}_j*f, \phi\rangle,
\quad \phi\in\cS_+.
$$
We now employ (\ref{rep-Psi*Phi-1}) and the Cauchy-Schwarz inequality
to obtain, for $j\ge 2$,
\begin{align*}
|\langle \Psi_j\ast\overline{\Phi}_j\ast f, \phi\rangle|^2
&=\Big|\sum_{m=4^{j-2}+1}^{4^j}
\overline{\ha\Big(\frac{m}{4^{j-1}}\Big)}\hb\Big(\frac{m}{4^{j-1}}\Big)
\langle \cF^\a_m*f, \cF^\a_m*\phi\rangle\Big|^2\\
&\le\sum_{m=4^{j-2}+1}^{4^j}
\Big|\ha\Big(\frac{m}{4^{j-1}}\Big)\Big|^2\|\cF^\a_m*f\|_2^2
\sum_{m=4^{j-2}+1}^{4^j}
\Big|\hb\Big(\frac{m}{4^{j-1}}\Big)\Big|^2\|\cF^\a_m*\phi\|_2^2\\
&\le \|\Phi_j*f\|_2^2
\sum_{m=4^{j-2}+1}^{4^j}
\|\cF^\a_m*\phi\|_2^2.
\end{align*}
Using inequality (\ref{norms3}) we get
\begin{align*}
\|\Phi_j*f\|_2
\le c2^{j(M+|s|)}\|2^{sj}\WW(2^j;\cdot)^{-\r/d}\Phi_j*f(\cdot)\|_p
\le c2^{j(M+|s|)}\|f\|_{\Fsrpq},
\end{align*}
where $M$ depends on $p$, $\alpha$, $d$, and $\rho$.
From the above estimates we infer
\begin{align*}
|\langle \Psi_j\ast\overline{\Phi}_j\ast f, \phi\rangle|
\le c 2^{-j}\|f\|_{\Fsrpq} 2^{jk} \sum_{4^{j-2}<m\le 4^j}\|\cF^\a_m*f\|_2
\le c 2^{-j}\|f\|_{\Fsrpq}P^*_k(\phi)
\end{align*}
for $k\ge M+|s|+1$.
A similar estimate trivially holds for $j=0, 1$.
Summing up we get
$$
|\langle f, \phi\rangle| \le c\|f\|_{\Fsrpq}P^*_k(\phi),
$$
which completes the proof.
$\qed$


\begin{prop}\label{prop:identification}
The following identification holds:
\begin{equation}\label{ident1}
F^{0 0}_{p 2} \sim L^p(\w),
\quad 1 < p < \infty,
\end{equation}
with equivalent norms.
\end{prop}
The proof of this proposition is the same as the proof of Proposition~4.3 in \cite{NPW}
in the case of spherical harmonics. We omit it.
Almost arbitrary $L^p$ multipliers for Laguerre expansions can be used for the proof.
However, since we cannot find in the literature any multipliers
for the Laguerre expansions we use in the present paper,
we next give easy to prove but non-optimal multipliers.


\begin{prop}\label{prop:multipliers}
Let $k$ be sufficiently large integer $(k > (5/2)|\a|+ (7/4)d+3$ will do$)$
and suppose $m \in C^k(\RR_+)$ obeys
\begin{equation}\label{cond-on-m}
\sup_{t\in\RR_+} |t^j m^{(j)}(t)| \le c
\quad \mbox{for \,$j=0, 1, \dots, k.$}
\end{equation}
Then the operator
$
T_m^\a(f):= \sum_{j=0}^\infty m(j) \cF_j^\a*f
$
is bounded on $L^p(\w)$, $1<p<\infty$.
\end{prop}

The proof is given in \S\ref{Proofs2}.

\subsection{Needlet Decomposition of Laguerre-Triebel-Lizorkin Spaces}\label{decomposition-F}

As a companion to $\Fsrpq$ we now introduce the sequence spaces $\fsrpq$.
Here  $\{\cX_j\}_{j=0}^\infty$ is the sequence of points from (\ref{def-Xj})
with associated tiles $\{R_\xi\}_{\xi\in\cX_j}$, defined in ({\ref{def-Rxi}).
Just as in the definition of needlets in \S\ref{def-needlets},
we set $\cX:=\cup_{j\ge 0} \cX_j$.


\begin{defn}
Suppose $s, \r \in \RR$, $0<p<\infty$, and $0<q\le\infty$. Then~$\fsrpq$
is defined as the space of all complex-valued sequences
$h:=\{h_{\xi}\}_{\xi\in \cX}$ such that
\begin{equation}\label{def-f-space}
\|h\|_{\fsrpq} :=\Big\|\Big(\sum_{j=0}^\infty
2^{sjq}\sum_{\xi \in \cX_j}
[|h_{\xi}|\WW(4^j;\xi)^{-\r/d}\tONE_{R_\xi}(\cdot)]^q\Big)^{1/q}\Big\|_p <\infty
\end{equation}
with the usual modification for $q=\infty$. Recall that
$\tONE_{R_\xi}:=\mu(R_\xi)^{-1/2}\ONE_{R_\xi}$.
\end{defn}

In analogy to the classical case on $\RR^d$ we introduce
``analysis" and ``synthesis" operators by
\begin{equation}\label{anal_synth_oprts}
S_\varphi: f\rightarrow \{\langle f, \varphi_\xi\rangle \}_{\xi \in \cX}
\quad\text{and}\quad
T_\psi: \{h_\xi\}_{\xi \in \cX}\rightarrow \sum_{\xi\in \cX}h_\xi\psi_\xi.
\end{equation}
We next show that the operator $T_\psi$ is well defined on $\fsrpq$.


\begin{lem}\label{lem:synthesis}
Let $s, \rho \in \RR$, $0<p<\infty$, and $0<q\le \infty$. 
Then for any $h\in  \fsrpq$,
$T_\psi h:=\sum_{\xi\in \cX}h_\xi \psi_\xi$ converges in $\cS_+'$.
Moreover, the operator
$T_\psi: \fsrpq \to \cS_+'$ is continuous, i.e.
there exist constants $k>0$ and $c>0$ such that
\begin{equation}\label{synthesis}
|\langle T_\psi h, \phi\rangle| \le c P^*_k(\phi)\|h\|_{\fsrpq}
\quad \mbox{for} \;\; h\in\fsrpq \;\;\mbox{and}\;\; \phi\in \cS_+.
\end{equation}
\end{lem}

\noindent
{\bf Proof.}
Let $h\in \fsrpq$. Using the definition of $\fsrpq$ we obtain
$$
2^{js} |h_\xi| \WW(4^j; \xi)^{-\rho/d}\|\tONE_{R_\xi}(\cdot)\|_p \le \|h\|_{\fsrpq}
\quad \mbox{for $\xi\in \cX_j$, \;$j\ge 0$.}
$$
But (\ref{size-R-xi3}) gives
$\|\tONE_{R_\xi}\|_p = \mu(R_\xi)^{1/p-1/2}
\ge c[2^{-jd} \WW(4^j, \xi)]^{1/p-1/2}$
for $\xi\in\cX_j$
and since $2^{-j(2|\alpha|+d)} \le \WW(4^j, \xi) \le c2^{j(2|\alpha|+d)}$
it follows that for $\xi\in\cX_j$
\begin{equation}\label{estim-h-xi}
|h_\xi|
\le c2^{jM} \|h\|_{\fsrpq}
\quad\mbox{with}\quad M:= |s|+2(|\alpha|+d)(|\rho|/d+|1/p-1/2|).
\end{equation}

By Lemma~\ref{lem:char-S}
$\phi=\sum_{n=0}^\infty \cF^\a_n*\phi$ in $\cS_+$ for $\phi\in\cS_+$
and hence
for $\xi\in\cX_j$
$$
\psi_\xi(x):=c_\xi^{1/2} \Psi_j(x, \xi)
=c_\xi^{1/2} \sum_{4^{j-2}<m<4^j}\hb\Big(\frac{m}{4^{j-1}}\Big)\PP_m(x, \xi),
\quad c_\xi\sim |R_\xi|\WW(4^j, \xi).
$$
Therefore,
$$
\langle \psi_\xi, \phi\rangle=c_\xi^{1/2}
\sum_{4^{j-2}<m<4^j}\hb\Big(\frac{m}{4^{j-1}}\Big)\cF^\a_m*\overline{\phi}
$$
and hence
$$
|\langle \psi_\xi, \phi\rangle|
\le c 2^{-j(|\alpha|+d)}\sum_{4^{j-2}<m<4^j}\|\cF^\a_m*\phi\|_\infty.
$$
Since $\cF^\a_m*\phi\in V_m$, by Proposition~\ref{prop:norms}
$\|\cF^\a_m*\phi\|_\infty \le cm^{(d+|\a|)/2)} \|\cF^\a_m*\phi\|_2$
and hence
$$
|\langle \psi_\xi, \phi\rangle|
\le c 2^{j(2|\alpha|+2d)}
\sum_{4^{j-2}<m<4^j}\|\cF^\a_m*\phi\|_2.
$$
This along with (\ref{estim-h-xi}) and the fact that
$\#\cX_j\le c4^{jd}$ yields, for $\phi\in\cS_+$,
\begin{align}\label{T-psi-bound}
\sum_{\xi\in\cX}|h_\xi| |\langle \psi_\xi, \phi\rangle|
&\le \sum_{j=0}^\infty \sum_{\xi\in\cX_j}|h_\xi||\langle \psi_\xi, \phi\rangle| \notag\\
&\le c \|h\|_{\fsrpq}\sum_{j=0}^\infty (\#\cX_j)2^{j(M+2|\alpha|+2d)}
\sum_{4^{j-2}<m<4^j}\|\cF^\a_m*\phi\|_2\\
&\le c \|h\|_{\fsrpq}\sum_{m=0}^\infty (m+1)^{k}\|\cF^\a_m*\phi\|_2
\sum_{j=0}^\infty 2^{j(M+2|\alpha|+4d+1-k)} \notag\\
&\le c \|h\|_{\fsrpq}P^*_k(\phi), \notag
\end{align}
where $k:= \lfloor M+2|\alpha|+4d+2\rfloor > M+2|\alpha|+4d+1$.
Therefore, the series
$\sum_{\xi\in\cX} h_\xi \psi_\xi$ converges in $\cS'$.
We define $T_\psi h$ by
$\langle T_\psi h, \phi \rangle := \sum_{\xi\in\cX} h_\xi\langle  \psi_\xi, \phi\rangle$
for all $\phi\in \cS$.
Estimate (\ref{synthesis}) follows by (\ref{T-psi-bound}).
$\qed$

\smallskip

We now present our main result on Laguerre-Triebel-Lizorkin spaces.


\begin{thm}\label{thm:needlet-Tr-Liz}
Let $s, \r\in \RR$, $0< p< \infty$ and $0<q\le \infty$.
Then the operators
$S_\varphi:\Fsrpq\rightarrow\fsrpq$ and $T_\psi:\fsrpq\rightarrow \Fsrpq$
are bounded and $T_\psi\circ S_\varphi=Id$ on $\Fsrpq$.
Consequently, $f\in \Fsrpq$ if and only if
$\{\langle f, \varphi_\xi\rangle \}_{\xi \in \cX}\in \fsrpq$
and
\begin{align}\label{Fnorm-equivalence-1}
\|f\|_{\Fsrpq}
\sim  \|\{\langle f,\varphi_\xi\rangle \}\|_{\fsrpq}.
\end{align}
In addition, the definition of $\Fsrpq$ is independent of the particular selection of
$\ha$ satisfying $(\ref{ha1})$--$(\ref{ha2})$.
\end{thm}

To prove this theorem we need several lemmas with proofs given
in Section \ref{Proofs2}.
Assume that $\{\Phi_j\}$ are the kernels from the definition
of Laguerre-Triebel-Lizorkin spaces
and $\{\ph_\xi\}_{\xi\in\cX}$ and $\{\psi_\xi\}_{\xi\in\cX}$
are needlet systems defined as in (\ref{def-needlets1})
with no connection between the functions $\ha$'s from (\ref{def-Phi-j}) and (\ref{def.Phi-j}).
We also assume that $p,q$ from the hypothesis of Theorem~\ref{thm:needlet-Tr-Liz}
are fixed and we choose $0<t<\min\{p, q\}$.


\begin{lem}\label{lem:Phi*psi}
For any $\sigma > d$ there exists a constant $c_\sigma>0$ such that
\begin{equation}\label{Phi*psi}
|\Phi_j\ast \psi_\xi(x)|
\le \frac{c_\sigma}{\mu(R_\xi)^{1/2}(1+2^m \|x-\xi\|)^{\sigma}},
\quad \xi\in\cX_m, \quad  j-1\le m \le j+1,
\end{equation}
and $\Phi_j\ast \psi_\xi\equiv 0$ for $\xi\in\cX_m$ if $|m-j|\ge 2$,
were $\cX_m:=\emptyset$ if $m < 0$.
\end{lem}


\begin{defn}\label{def-s*}
For any collection of complex numbers $\{h_\xi\}_{\xi \in \cX_j}$ $(j\ge 0)$,
we define
\begin{equation}\label{def-s}
h^*_j(x) := \sum_{\eta \in \cX_j} \frac{|h_\eta|}{(1+ 2^{j}\|\eta-x\|)^\lambda}
\end{equation}
and
\begin{equation}\label{def-a*}
h_\xi^*:=h^*_j(\xi),
\quad \xi \in \cX_j,
\end{equation}
where $\lambda := 2d+2(|\alpha|+3d)/t+2(|\alpha|+d)|\rho|/d$.
\end{defn}


\begin{lem}\label{lem:sum<M}
For any set $\{h_\eta\}_{\eta \in \cX_j}$ $(j \ge 0)$ of complex numbers
\begin{equation}\label{s*<M}
h_j^*(x)
\le c\cM_t\Big(\sum_{\eta \in \cX_j}|h_\eta| \ONE_{R_\eta}\Big)(x),
\quad x\in\RR^d_+.
\end{equation}
Moreover, for $\xi \in \cX_j$
\begin{equation}\label{sum<M}
\WW(4^j;\xi)^{-\rho/d} h_\xi^* \ONE_{R_\xi}(x)
\le c\cM_t\Big(\sum_{\eta \in \cX_j}|h_\eta| \WW(4^j;\eta)^{-\rho/d} \ONE_{R_\eta}\Big)(x),
\;\;x\in\RR^d_+.
\end{equation}
Here the constants depend only on $d$, $\alpha$, $\rho$, $\delta$, and $t$.
\end{lem}


\begin{lem}\label{lem:a*=b*}
Suppose $g \in \V_{4^j}$ and denote
$$
M_\xi:=\sup_{x \in R_\xi} |g(x)|,
\quad \xi\in\cX_j,
\quad\mbox{and}\quad
m_\eta:=\inf_{x \in R_\eta} |g(x)|, \quad
\eta \in \cX_{j+\ell}.
$$
Then there exists $\ell \ge 1$, depending only $d$, $\alpha$, $\delta$, and $\lambda$, such that
for any $\xi \in \cX_j$
\begin{equation}\label{a*=b*}
M_\xi^* \le c m_\eta^*
\quad \mbox{for all \; $\eta\in\cX_{j+\ell}$, $R_\eta\cap R_\xi\ne \emptyset$},
\end{equation}
and, therefore,
\begin{equation}\label{M*<m*2}
M_\xi^*\ONE_{R_\xi}(x)
\le c \sum_{\eta\in\cX_{j+\ell}, R_\eta\cap R_\xi\ne \emptyset}
m_\eta^*\ONE_{R_\eta}(x),
\quad x\in \R^d,
\end{equation}
where $c>0$ depends only on $d$, $\alpha$, $\delta$, and $t$.
\end{lem}

\noindent
{\bf Proof of Theorem \ref{thm:needlet-Tr-Liz}.}
Choose $\sigma$ so that $\sigma \ge \lambda+2(|\alpha|+d)|\rho|/d$
and recall that $t$ has already been selected so that $0 < t < \min\{p, q\}$.

Suppose $\{\Phi_j\}$ are from the definition
of Laguerre-Triebel-Lizorkin spaces (see (\ref{def-Phi-j})-(\ref{ha2})).
As already mentioned in \S\ref{Def-of-needlets},
there exists a function $\hb$ satisfying (\ref{ab1})-(\ref{ab3})
such that (\ref{ab5}) holds as well.
Using this function we define $\{\Psi_j\}$ just as in (\ref{def-Psi-j}).
Then we use $\{\Phi_j\}$ and $\{\Psi_j\}$ to define as in (\ref{def-needlets1})
a~pair of dual needlet systems $\{\ph_\eta\}$ and $\{\psi_\eta\}$.

Suppose $\{\wtd\ph_\eta\}$, $\{\wtd\psi_\eta\}$ is a second pair of needlet systems,
defined as in (\ref{def.Phi-j})-(\ref{def-needlets1}) using another pair of
kernels $\{\wtd\Phi_j\}$, $\{\wtd\Psi_j\}$.

We first show the boundedness of the operator
$T_{\wtd\psi}: \fsrpq \to \Fsrpq$.
Let $h\in\fsrpq$ and set $f:=T_{\wtd\psi} h= \sum_{\xi\in \cX} h_\xi\wtd\psi_\xi$.
Evidently $\Phi_j*\wtd\psi_\xi =0$ if $\xi\in\cX_m$
and $|j-m|\ge 2$, and hence
$$
\Phi_j*f = \sum_{m=j-1}^{j+1}\sum_{\xi \in \cX_m}
h_\xi \Phi_j*\wtd\psi_\xi
\qquad (\cX_{-1}:=\emptyset).
$$
Denote $H_\xi:= h_\xi \WW(4^m;\xi)^{-\rho/d}\mu(R_\xi)^{-1/2}$. 
Using Lemma~\ref{lem:Phi*psi} and (\ref{W<W}) we get
\begin{align}\label{Phij-f<h*}
\WW(4^j;x)^{-\rho/d}|\Phi_j*f(x)|
&\le \sum_{m=j-1}^{j+1}\sum_{\xi \in \cX_m}|h_\xi|\WW(4^j;x)^{-\rho/d}|\Phi_j*\wtd\psi_\xi(x)|\notag\\
&\le c\sum_{m=j-1}^{j+1}\sum_{\xi \in \cX_m}
\frac{|h_\xi|\WW(4^m;\xi)^{-\rho/d}\mu(R_\xi)^{-1/2}}
{(1+2^m \|\xi-x\|)^{\sigma-2(|\alpha|+d)|\rho|/d}}\\
&\le c\sum_{m=j-1}^{j+1}H_m^*(x)
\qquad (H_{-1}^*:=0), \notag
\end{align}
where
$H_m^*(x)$ is defined as in (\ref{def-s}).
We use this in the definition of $\|f\|_{\Fsrpq}$ and apply
Lemma \ref{lem:sum<M} and the maximal inequality (\ref{max-ineq}) to obtain
\begin{align*}
\|f\|_{\Fsrpq}
&\le \Big\|\Big(\sum_{j=0}^\infty
(2^{js}|H_j^*(\cdot)|)^q\Big)^{1/q}\Big\|_p\notag\\
&\le c\Big\|\Big(\sum_{j=0}^\infty
\Big[\cM_t\Big(2^{js}\sum_{\xi\in\cX_j}
|h_\xi|\WW(4^j;\xi)^{-\rho/d}\mu(R_\xi)^{-1/2}\ONE_{R_\xi}\Big)
\Big]^q\Big)^{1/q}\Big\|_p\\
&\le c\|\{h_\eta\}\|_{\fsrpq}.
\end{align*}
Hence the operator $T_{\wtd\psi}: \fsrpq\to\Fsrpq$
is bounded.


Let the space $\Fsrpq$ be defined using $\{\overline{\Phi}_j\}$
instead of $\{\Phi_j\}$. We now prove the boundedness of the
operator $S_\ph: \Fsrpq \to \fsrpq$.
Let $f\in \Fsrpq$ and denote
$$
M_\xi := \sup_{x \in R_\xi}|\overline{\Phi}_j*f(x)|,
\;\; \xi \in \cX_j,
\quad\mbox{and}\quad
m_\eta :=\inf_{x \in R_\eta}|\overline{\Phi}_j*f(x)|,
\;\; \eta \in \cX_{j+\ell},
$$
where $\ell$ is the constant from Lemma~\ref{lem:a*=b*}.
We have
\begin{equation}\label{coef<M*}
|\langle f, \ph_\xi\rangle|
\le c_\xi^{1/2}|\overline{\Phi}_j*f(\xi)|
\le c\mu(R_\xi)^{1/2}M_\xi
\le c\mu(R_\xi)^{1/2} M_\xi^*.
\end{equation}
Evidently, $\overline{\Phi}_j*f \in \V_{4^j}$,
and applying Lemma~\ref{lem:a*=b*} (see (\ref{M*<m*2})), we get
\begin{equation}\label{M*<m*}
M_\xi^*\ONE_{R_\xi}(x)
\le c \sum_{\eta\in\cX_{j+\ell}, R_\eta\cap R_\xi\ne \emptyset}
m_\eta^*\ONE_{R_\eta}(x),
\quad x\in \R^d.
\end{equation}
It is easy to see that
$\WW(4^{j+\ell}; y) \sim \WW(4^j; \xi)$ for $y\in R_\xi$.
We use this, (\ref{coef<M*})-(\ref{M*<m*}),
Lemma~\ref{lem:sum<M}, and the maximal inequality (\ref{max-ineq})
to obtain
\begin{align*}
\|\{\langle f, \ph_\xi \rangle\}\|_\fsrpq
&\le
c\Big\|\Big(\sum_{j=0}^\infty 2^{sjq}
\Big(\sum_{\xi \in \cX_{j}}
\WW(4^j; \xi)^{-\r/d}M_\xi^*\ONE_{R_\xi}\Big)^q\Big)^{1/q}\Big\|_p\\
&\le
c\Big\|\Big(\sum_{j=0}^\infty 2^{sjq}
\Big(\sum_{\eta \in \cX_{j+\ell}}
\WW(4^{j+\ell}; \eta)^{-\r/d}m_\eta^*\ONE_{R_\eta}\Big)^q\Big)^{1/q}\Big\|_p\\
&\le c\Big\|\Big(\sum_{j=0}^\infty
\cM_t\Big(2^{sj}\sum_{\eta \in \cX_{j+\ell}}
\WW(4^{j+\ell}; \eta)^{-\r/d} m_\eta\ONE_{R_\eta}\Big)^q\Big)^{1/q}\Big\|_p\\
&\le c\Big\|\Big(\sum_{j=0}^\infty
\Big(2^{sj}\sum_{\eta \in \cX_{j+\ell}}
\WW(4^{j+\ell}; \eta)^{-\r/d} m_\eta\ONE_{R_\eta}\Big)^q\Big)^{1/q}\Big\|_p\\
&\le c\Big\|\Big(\sum_{j=0}^\infty 2^{sjq}
\WW(4^j; \cdot)^{-\r/d}|\overline{\Phi}_j*f(\cdot)|^q \Big)^{1/q}\Big\|_p
= c\|f\|_\Fsrpq.
\end{align*}
Here for the second inequality we used that each tile $R_\eta$, $\eta\in\cX_{j+l}$,
intersects no more that finitely many (depending only on $d$) tiles $R_\eta$, $\eta\in\cX_j$.
The above estimates prove the boundedness of the operator
$S_\ph: \Fsrpq \to \fsrpq$.
The identity $T_\psi\circ S_\ph=Id$ follows by Theorem~\ref{prop:needlet-rep}.

It remains to show the independence of the definition of Triebel-Lizorkin spaces
from the specific selection of $\ha$ satisfying (\ref{ha1})-(\ref{ha2}).
Suppose $\{\Phi_j\}$, $\{\wtd\Phi_j\}$ are two sequences of
kernels as in the definition of Triebel-Lizorkin spaces defined by
two different functions $\ha$ satisfying (\ref{ha1})-(\ref{ha2}).
As above there exist two associated needlet systems
$\{\Phi_j\}$, $\{\Psi_j\}$, $\{\ph_\xi\}$, $\{\psi_\xi\}$
and
$\{\wtd\Phi_j\}$, $\{\wtd\Psi_j\}$, $\{\wtd\ph_\xi\}$, $\{\wtd\psi_\xi\}$.
Denote by $\|f\|_{\Fsrpq(\Phi)}$ and
$\|f\|_{\Fsrpq(\wtd\Phi)}$ the $F$-norms defined via
$\{\Phi_j\}$ and $\{\wtd\Phi_j\}$.
Then from above
$$
\|f\|_{\Fsrpq(\Phi)}
\le c \|\{\langle f, \wtd\ph_\xi\rangle\}\|_{\fsrpq}
\le c\|f\|_{\Fsrpq(\overline{\wtd\Phi})}.
$$
The independence of the definition of $\Fsrpq$ of the specific
choice of $\ha$ in the definition of the functions $\{\Phi_j\}$
follows by interchanging the roles of $\{\Phi_j\}$ and $\{\wtd\Phi_j\}$
and their complex conjugates.
$\qed$

\smallskip

To us the spaces $\F sspq$ are more natural than the spaces
$\Fsrpq$ with $\r\ne s$ since they  embed correctly with respect
to the smoothness index $s$.

%
\begin{prop}\label{prop:F-embedding}
Let $0<p<p_1<\infty$, $0<q, q_1\le\infty$, and
$-\infty<s_1<s<\infty$. Then we have the continuous embedding
\begin{equation}\label{F-embed}
\F sspq \subset \F {s_1}{s_1}{p_1}{q_1} \quad\mbox{if}\quad
s/d-1/p=s_1/d-1/p_1.
\end{equation}
\end{prop}

The proof of this embedding result can be carried out similarly
as the proof of Proposition~4.11 in \cite{KPX2},
using the idea of the proof in the classical case on $\RR^n$
(see e.g. \cite{T1}, page 129). We omit it.

\section{Laguerre-Besov Spaces}\label{Besov}
\setcounter{equation}{0}

We introduce weighted Besov spaces on $\RR^d_+$ in the context of
Laguerre expansions using the kernels $\{\Phi_j\}$ from (\ref{def-Phi-j}) with
$\ha$ satisfying (\ref{ha1})-(\ref{ha2})
(see \cite{Peetre}, \cite{T1} for the general idea of using orthogonal or
spectral decompositions in defining Besov spaces).

\subsection{Definition of Laguerre-Besov Spaces}


\begin{defn}\label{def-Besov}
Let $s, \r\in \RR$ and $0<p,q \le \infty$. The Laguerre-Besov space
$\Bsrpq := \Bsrpq(\cF^\a)$ is defined
as the set of all $f \in \cS_+'$ such that
\begin{equation}\label{def-Besov-sp}
\|f\|_{\Bsrpq} :=
\Big(\sum_{j=0}^\infty \Big(2^{s j}
\|\WW(4^j; \cdot)^{-\r/d}\Phi_j*f(\cdot)\|_p\Big)^q\Big)^{1/q}
< \infty,
\end{equation}
where the $\ell_q$-norm is replaced by the sup-norm if $q=\infty$.
\end{defn}

Observe that as in the case of Laguerre-Triebel-Lizorkin spaces the above definition
is independent of the particular choice of $\ha$ obeying  (\ref{ha1})-(\ref{ha2})
(see Theorem~\ref{thm:Bnorm-eq}).
Also, as for $\Fsrpq$ the Besov space $\Bsrpq$ is a quasi-Banach
space which is continuously embedded in $\cS_+'$. We skip the details.

\subsection{Needlet Decomposition of Laguerre-Besov Spaces}

We next define the sequence spaces $\bsrpq$
associated to the Laguerre-Besov spaces $\Bsrpq$.
As in \S\ref{Tri-Liz} we assume that $\{\cX_j\}_{j=0}^\infty$ are from (\ref{def-Xj})
with associated tiles $\{R_\xi\}_{\xi\in\cX_j}$ from (\ref{def-Rxi}).
As before we set $\cX:=\cup_{j\ge 0} \cX_j$.


\begin{defn}
Let $s, \r\in \RR$ and $0<p,q \le \infty$. Then $\bsrpq$
is defined to be the space of all complex-valued sequences
$h:=\{h_{\xi}\}_{\xi\in \cX}$ such that
\begin{equation}\label{def-berpq}
\|h\|_{\bsrpq} :=\Bigl(\sum_{j=0}^\infty
2^{jsq}
\Bigl[\sum_{\xi\in \cX_j}\Big(\WW(4^j;\xi)^{-\r/d}\mu(R_\xi)^{1/p-1/2}|h_\xi|\Big)^p
\Bigr]^{q/p}\Bigr)^{1/q}
\end{equation}
is finite,
with the usual modification whenever $p=\infty$ or $q=\infty$.
\end{defn}

We shall utilize again the analysis and synthesis operators
 $S_\ph$ and $T_\psi$ defined in (\ref{anal_synth_oprts}).
The next lemma guarantees that the operator $T_\psi$ is well defined on $\bsrpq$.


\begin{lem}\label{lem:synthesis-B}
Let $s, \rho \in \RR$ and $0<p, q\le\infty$.
Then for any $h\in  \bsrpq$,
$T_\psi h:=\sum_{\xi\in \cX}h_\xi \psi_\xi$ converges in $\cS_+'$.
Moreover, the operator
$T_\psi: \bsrpq \to \cS_+'$ is continuous.
\end{lem}
The proof of this lemma is quite similar to the proof of Lemma~\ref{lem:synthesis}
and will be omitted.

Our main result in this section is the following characterization of
Laguerre-Besov spaces.


\begin{thm}\label{thm:Bnorm-eq}
Let $s, \r\in \RR$ and  $0< p,q\le \infty$.
Then the operators
$S_\varphi:\Bsrpq\rightarrow\bsrpq$ and
$T_\psi:\bsrpq\rightarrow \Bsrpq$
are bounded and $T_\psi\circ S_\varphi=Id$ on $\Bsrpq$.
Consequently, for $f\in \cS_+'$ we have that $f\in \Bsrpq$ if and
only if $\{\langle f, \varphi_\xi\rangle\}_{\xi \in \cX}\in \bsrpq$
and
\begin{align}\label{Bnorm-equivalence-1}
\|f\|_{\Bsrpq} & \sim  \|\{\langle f,\varphi_\xi\rangle\}\|_{\bsrpq}.
\end{align}
In addition, the definition of $\Bsrpq$ is independent of the particular
selection of $\ha$ satisfying $(\ref{ha1})$--$(\ref{ha2})$.
\end{thm}

The proof of this theorem relies on some lemmas from the proof of
Theorem~\ref{thm:needlet-Tr-Liz} as well as the next lemma
with proof given in Section~\ref{Proofs2}.


\begin{lem}\label{l:half_shannon}
Let $0<p\le \infty$ and $\rho\in \RR$.
Then for any $g\in V_{4^j}, j\ge 0$,
\begin{equation}
 \Big(\sum_{\xi\in \cX_j}\WW(4^j; \xi)^{-\rho p/d}\max_{x\in R_\xi}|g(x)|^p \mu(R_\xi)\Big)^{1/p}
 \le c\|\WW(4^j; \cdot)^{-\rho/d} g(\cdot)\|_p.
\end{equation}
\end{lem}


\noindent
{\bf Proof of Theorem~\ref{thm:Bnorm-eq}.}
We will use some basic assumptions and notation from the proof
of Theorem~\ref{thm:needlet-Tr-Liz}.
Let $0<t<p$ and $\sigma \ge \lambda+2(|\alpha|+d)|\rho|/d$.
Assume that
$\{\Phi_j\}$, $\{\Psi_j\}$, $\{\ph_\eta\}$, $\{\psi_\eta\}$
and
$\{\wtd\Phi_j\}$, $\{\wtd\Psi_j\}$, $\{\wtd\ph_\eta\}$, $\{\wtd\psi_\eta\}$
are two needlet systems, defined as in (\ref{def.Phi-j})-(\ref{def-needlets1}),
that originate from two completely different functions $\ha$ satisfying
(\ref{ha1})-(\ref{ha2}).

Let us first prove the boundedness of the operator
$T_{\wtd\psi}: \bsrpq \to \Bsrpq$,
assuming that $\Bsrpq$ is defined by $\{\Phi_j\}$.
Suppose $h\in \bsrpq$ and set
$f:= T_{\wtd\psi} h = \sum_{\xi\in \cX} h_\xi\wtd\psi_\xi.$

Denote $H_\xi:= h_\xi \WW(4^m;\xi)^{-\rho/d}\mu(R_\xi)^{-1/2}$, $\xi\in\cX_m$.
Then by (\ref{Phij-f<h*})
and Lemma~\ref{lem:sum<M}
\begin{align*}
&\|\WW(4^j;\cdot)^{-\rho/d}\Phi_j*f(\cdot)\|_p
\le c \sum_{m=j-1}^{j+1}\|H_m^*\|_p\\
& \qquad\le c \sum_{m=j-1}^{j+1}\Big\|\cM_t\Big(\sum_{\xi\in\cX_m}
|h_\xi|\WW(4^m;\xi)^{-\rho/d}\mu(R_\xi)^{-1/2}\ONE_{R_\xi}\Big)\Big\|_p\\
& \qquad\le c \sum_{m=j-1}^{j+1}\Big(\sum_{\xi\in\cX_m}
\Big(|h_\xi|\WW(4^m;\xi)^{-\rho/d}\mu(R_\xi)^{1/p-1/2}\Big)^p\Big)^{1/p},
\end{align*}
which yields
$\|f\|_{\Bsrpq} \le c \|\{h_\eta\}\|_{\bsrpq}$
and hence the claimed boundedness of $T_{\wtd\psi}$.


We now prove the boundedness of the operator
$S_\varphi:\Bsrpq\rightarrow \bsrpq$,
where we assume that the space $\Bsrpq$ is defined in terms
of $\{\overline{\Phi}_j\}$ in place of $\{\Phi_j\}$.
Just as in (\ref{coef<M*}) we have
$|\langle f,\varphi_\xi\rangle|
\le c\mu(R_\xi)^{1/2}|\overline{\Phi}_j*f(\xi)|$,
$\xi \in \cX_j$.
Since $\overline{\Phi}_j\ast f\in V_{4^j}$, Lemma~\ref{l:half_shannon} implies
\begin{align*}
&\sum_{\xi\in \cX_j}\Big(\WW(4^j;\xi)^{-\r/d} \mu(R_\xi)^{1/p-1/2}|\langle f,\varphi_\xi\rangle|\Big)^p\\
&\qquad \le c\sum_{\xi\in \cX_j}\WW(4^j;\xi)^{-\r p/d}|\overline{\Phi}_j*f(\xi)|^p\mu(R_\xi)
\le c\|\WW(4^j;\cdot)^{-\r/d}\overline{\Phi}_j*f(\cdot)\|_p^p,
\end{align*}
which leads immediately to
$\|\{\langle f, \ph\rangle\}\|_{\bsrpq} \le c\|f\|_{\Bsrpq}$.

The identity $T_\psi\circ S_\varphi=Id$ follows by Proposition~\ref{prop:needlet-rep}.
The independence of $\Bsrpq$ of the specific selection of $\ha$
in the definition of  $\{\Phi_j\}$ follows from above exactly as
in the Triebel-Lizorkin case (see the proof of
Theorem~\ref{thm:needlet-Tr-Liz}).
$\qed$

\smallskip

The parameter $\rho$ in the definition of the Besov spaces $\Bsrpq$ allows one to
consider various scales of spaces.
A ``classical" choice of $\r$ would be $\r=0$.
However, to us most natural are the spaces $\B sspq$ ($\r=s$)
for they embed ``correctly" with respect to
the smoothness index $s$:

%
\begin{prop}\label{B-embedding}
Let $0<p\le p_1\le\infty$, $0<q\le q_1\le \infty$, and
$-\infty<s_1\le s<\infty$. Then we have the continuous embedding
\begin{equation}\label{B-embed}
\B sspq \subset \B {s_1}{s_1}{p_1}{q_1} \quad\mbox{if}\quad
s/d-1/p=s_1/d-1/p_1.
\end{equation}
\end{prop}

\noindent
{\bf Proof.}
Assuming that $\Phi_j$ is from Definition~\ref{def-Besov} we have
$\Phi_j*f\in V_{4^{j+1}}$ and applying estimate (\ref{norms2}) from
Proposition~\ref{prop:norms} we obtain
$$
\|\WW (4^j;\cdot)^{-s_1/d}\Phi_j*f(\cdot)\|_{p_1}
\le c2^{jd(1/p-1/p_1)}\|\WW(4^j;\cdot)^{-s/d}\Phi_j*f(\cdot)\|_p,
$$
where we used that $s/d-1/p=s_1/d-1/p_1$.
This implies (\ref{B-embed}) at once.
$\qed$


\section{Proofs for Sections \ref{localization}-\ref{background}}\label{Proofs}
\setcounter{equation}{0}

\subsection{Proof of estimates (\ref{est-Lambda-n}) and (\ref{est-Lambda-n-large})
in Theorem \ref{thm:est-Lambda-n}}
We may assume that $n\ge n_0$, where $n_0$ is a sufficiently large constant.
Estimate (\ref{est-Lambda-n}) will be established by applying repeatedly summation
by parts to the sum in the definition (\ref{def.Ln}) of $\LL_n(x,y)$. For
a sequence of numbers $\{a_m\}$ we denote by $\Delta^k a_m$ the $k$th
forward differences, defined by $\Delta a_m:=a_m-a_{m+1}$ and inductively
$\Delta^{k+1} a_m:= \Delta(\Delta^{k}a_m)$.
Choose $k\ge \sigma+2|\a|+d/2$ and denote
\begin{equation} \label{def.Sigma-n}
\Omega_n^k(x, y):=\sum_{m=0}^n A_{n-m}^k\PsiP_m^\a(x, y),
\quad A_m^k:=\binom{m+k}{k}.
\end{equation}
Using summation by parts $k$ times, we obtain
\begin{equation} \label{rep1-Ln}
\LL_n(x, y)
:= \sum_{j=0}^\infty \ha\Big(\frac{j}{n}\Big) \PsiP_j^\a(x, y)
= \sum_{m=0}^\infty \Delta^{k+1}\ha\Big(\frac{m}{n}\Big)\cdot\Omega_m^k(x, y),
\end{equation}
where $\Delta^{k+1}$ is applied with respect to $m$.
By (\ref{cesaro-kernel}) and (\ref{def.Sigma-n}), it easily follows that
$\Omega_m^k (x,y) = c e^{- (\|x\|_2^2 + \|y\|_2^2)/2} P_m^{\a,k}(x^2, y^2)$
and combining this with (\ref{cesaroP}) we get
\begin{align*}
\Omega_m^k(x, y)&  = c \int_{[0,\pi]^d }
 L_m^{|\a|+k+d}\Big (\|x\|_2^2+\|y\|_2^2+2\sum_{i=1}^d x_iy_i\cos\theta_i \Big)\\
& \quad \times e^{-(\|x\|_2^2+\|y\|_2^2+2 \sum_{i=1}^d x_iy_i\cos \theta_i)/2}
 \prod_{i=1}^d j_{\a_i-1/2}(x_iy_i\cos\theta_i)\sin^{2\a_i}\theta_i \,d\theta.
\end{align*}
Using this in (\ref{rep1-Ln}) we arrive at the identity
\begin{align}\label{rep2-Ln}
\LL_n(x, y)
&= c  \int_{[0,\pi]^d } \KK_n^\lambda
     \Big (\|x\|^2+\|y\|^2+2\sum_{i=1}^d x_iy_i\cos\theta_i \Big)\\
& \qquad \times \prod_{i=1}^d j_{\a_i-1/2}(x_iy_i\cos\theta_i)
      \sin^{2\a_i}\theta_i \,d\theta, \notag
\end{align}
where $\lambda:=|\a|+k+d$ and the kernel $\KK_n^\lambda$ is defined by
\begin{equation} \label{Def-KK}
\KK_n^\lambda(t): = \sum_{m=0}^\infty \Delta^{k+1}\ha\Big(\frac{m}{n}\Big)
   L_m^{\lambda} (t) e^{- t/2}.
\end{equation}
By a well known property of finite differences we have
\begin{equation} \label{est-Delta-a}
\Big|\Delta^{k+1}\ha\Big(\frac{m}{n}\Big)\Big|
=n^{-k-1}|\ha^{(k+1)}(\xi)|
\le n^{-k-1}\|\ha^{(k+1)}\|_{L^\infty}.
\end{equation}
Further, it is known that \cite[p. 204]{AAR}
\begin{equation} \label{Bessel}
j_{\a-\frac12}(x)=x^{-\a+1/2}J_{\a-1/2}(x)
= \frac{2^{-\a+\frac12}}{\sqrt{\pi} \Gamma(\a)} \int_{-1}^1 e^{i x t}
        (1-t^2)^{\a -1} dt, \quad \a > 0,
\end{equation}
and $j_{-\frac12} (x) = \sqrt{\frac2\pi} \cos x$.
Therefore,
\begin{equation} \label{est-j-alpha}
        |j_{\a-\frac12}(x)|\le c_\a <\infty, \quad x\in \RR_+,\: \alpha \ge 0.
\end{equation}
By (\ref{est-Delta-a}) and (\ref{est3-eLn}) (with $\alpha$
replaced by $|\a|+k+1$) we obtain for $t > 0$
\begin{equation} \label{est-Fnk}
|\KK_n^\lambda(t) |
\le  c\sum_{m=\max\{\lceil un \rceil-k,1\}}^{\lfloor(1+v)n \rfloor}\frac{1}{m^{k+1}}
\left(\frac{m}{t}\right)^{(|\a|+k+d)/2} \le c n^{(-k+|\a|+d)/2} t^{-(|\a|+k+d)/2}.
\end{equation}
Using this in (\ref{rep2-Ln}) we get 
\begin{align*}
|\LL_n(x, y)|
  \le c n^{(-k+|\a|+d)/2}\int_{[0,\pi]^d}\frac{  \prod_{i=1}^d\sin^{2\a_i}\theta_i d\theta}
{\left(\|x\|_2^2+\|y\|_2^2+2 \sum_{i=1}^d x_iy_i\cos \theta_i\right)^{(k+|\a|+d)/2}}.
\end{align*}
Set $\tau:= (k+|\a|+d)/2$. Substituting $\theta_i = \pi-t_i$ in the above
integral and using $1- \cos t = 2 \sin^2 \frac t2 \sim t^2$ we infer
\begin{align} \label{eq:Jn}
|\LL_n(x, y)|
&  \le c n^{ -k+\tau}\int_{[0,\pi]^d}\frac{  \prod_{i=1}^d\sin^{2\a_i} t_i d t}
{\left(\|x - y\|_2^2+ 4 \sum_{i=1}^d x_iy_i\sin^2 \frac{t_i}2 \right)^{\tau}} \\
& \le  c n^{ -k+\tau} \int_{[0,\pi]^d}\frac{\prod_{i=1}^d t_i^{2\a_i} d t}
{\left(\|x - y\|^2+ \sum_{i=1}^d x_iy_i  t_i^2 \right)^{\tau}}
=:c M_n^{k, \a}(x, y). \notag
\end{align}
%
%
We estimate the integral above in two ways.
First, we trivially have
\begin{equation}\label{est2-Lambda-n}
|\LL_n(x, y)|\le c M_n^{k, \a}(x, y)
\le \frac{ cn^{-k + \tau}}{\|x-y\|^{2\tau}}
       \le \frac{cn^{|\a|+d}}{(n^{1/2}\|x-y\|)^{k+|\a|+d}}.
\end{equation}
%
%
The second estimate is really many estimates rolled into one.
For a fixed $1 \le \ell \le d$
we partition $\a$ into $\a=(\a', \a'')$ with
$\a' = (\a_1, \dots, \a_\ell)$ and $\a{''} = (\a_{\ell+1}, \dots, \a_d)$.
Since $\tau>|\alpha|+d/2$ and $x_iy_i>0$ we have
\begin{align*}
M_n^{k, \a}(x, y)
&\le cn^{-k+\tau}\int_{[0,\pi]^\ell}\frac{\prod_{i=1}^\ell t_i^{2\a_i} d t}
{\left(\|x - y\|^2+ \sum_{i=1}^\ell x_iy_i  t_i^2 \right)^{\tau}}\\
& \le \frac{cn^{-k+\tau}}{\prod_{i=1}^\ell (x_i y_i)^{\a_i+1/2}}  \prod_{i=1}^\ell
\int_0^{\pi(x_iy_i)^{1/2}}\frac{du}{(\|x-y\|^2+\sum_{i=1}^\ell u_i^2)^{\tau-|\a'|}},
\end{align*}
where we applied the substitutions $u_i=t_i(x_iy_i)^{1/2}$
and used $|\alpha'|$ power of the main term in the denominator to cancel the numerator.
Enlarging the integral domain to $\RR^\ell$ and using polar coordinates, the above
product of integrals is bounded by
\begin{align*}
\int_{\RR^\ell} \frac{du}{(\|x-y\|^2+\|u\|_2^2)^{\tau-|\a'|}}
=  \int_0^\infty \frac{r^{\ell-1} dr } {(\|x-y\|^2+r^2)^{\tau-|\a'|}}
\le   \frac{c}{\|x-y\|^{2(\tau - |\a'|)-\ell}}.
\end{align*}
From above and a little algebra we obtain for $1 \le \ell \le d$
\begin{align}\label{est1-Lambda-n}
&|\LL_n(x, y)|
\le cM_n^{k, \a}(x, y)\\
&\qquad\le \frac{c n^{d/2}} {\prod_{i=1}^\ell (x_iy_i)^{\a_i+\frac12}
\prod_{i=\ell+1}^d (n^{-1})^{\a_i+\frac12}(n^{1/2}\|x-y\|)^{k+|\a|-2|\a'|+d-\ell}}.\notag
\end{align}


A third bound on $|\LL_n(x,y)|$ will be obtained by estimating
all terms in (\ref{def.Ln}).
By (\ref{est2-eLn}) and (\ref{rel-Psi-L}) it follows that
\begin{equation}\label{Est-norm-F}
\|\PsiP_{\nu_i}^{\a_i}\|_\infty \le c \nu_i^{\a_i/2}, \; 1\le i\le d,
\quad\mbox{and}\quad
\|\PsiP_\nu^\a\|_\infty \le c \nu^{\a/2},
\end{equation}
and hence
$$
|\PsiP_m^\a(x,y)| \le c \sum_{|\nu| =m} \nu^{\a}
= c \binom{m+d-1}{m} m^{|\a|}\le c m^{|\a| + d -1},
\quad\mbox{yielding}
$$
\begin{equation}\label{est3-Lambda-n}
|\LL_n(x, y)| \le c\sum_{m=0}^{\lfloor(1+v)n\rfloor}m^{|\a|+d-1}  \le cn^{|\a|+d}.
\end{equation}


We also need the estimate:
\begin{equation}\label{est4-Lambda-n}
|\LL_n(x, y)| \le \frac{cn^{d/2} }
{ \prod_{i=1}^\ell (x_iy_i)^{\a_i+1/2} \prod_{i=\ell+1}^d (n^{-1})^{\a_i+1/2}},
\quad 1\le \ell\le d.
\end{equation}
By (\ref{est1-eLn}) it follows that
\begin{equation}\label{Est-Fna1}
|\cF_n^\a(x)|\le \frac{c}{x^{\a+1/2}n^{1/4}},
\quad \mbox{if \; $x^2\in \RR_+\setminus (2n+2\a+2, 6n+3\a+3)$},
\end{equation}
and if $x^2\in [2n+2\a+2, 6n+3\a+3)$ by (\ref{est4-eLn})
\begin{equation}\label{Est-Fna2}
|\cF_n^\a(x)|\le \frac{c}{x^{\a}n^{1/4}(n^{1/3}+|4n+2\alpha+2-x^2|)^{1/4}}.
\end{equation}
From these two estimates one easily concludes that for $x>0$
\begin{equation}\label{est-Fa-1}
|\cF_n^\a(x)|\le \frac{c}{x^{\a+1/2}n^{1/4}},
\quad \mbox{if \; $n\in \RR_+\setminus (x^2/5, x^2/3)$},
\end{equation}
and
\begin{equation}\label{est-Fa-2}
|\cF_n^\a(x)|\le \frac{c}{x^{\a+1/2}(1+|4n-x^2|)^{1/4}},
\quad \mbox{if \;$n\in [x^2/5, x^2/3]$.}
\end{equation}
Hence, $|\cF_n^\a(x)|$ can be bounded by the sum of the right-hand-side
quantities in (\ref{est-Fa-1})-(\ref{est-Fa-2}).
Also, from (\ref{Est-norm-F})
$\|\PsiP_{\nu_i}^{\a_i}\|_\infty \le c \nu_i^{\a_i/2}$.
From these along with (\ref{def.dPsi_n})-(\ref{def.dL_n})
we obtain
\begin{align*}
&|\LL_n(x, y)|
\le \prod_{i=1}^d \sum_{\nu_i=0}^{\lfloor(1+v)n\rfloor}
|\PsiP_{\nu_i}^{\a_i}(x_i)||\PsiP_{\nu_i}^{\a_i}(y_i)|
\le \frac{c}{\prod_{i=1}^\ell(x_i y_i)^{\a_i+1/2}}
\prod_{i=\ell+1}^d \sum_{\nu_i=0}^{\lfloor(1+v)n\rfloor}(\nu_i+1)^{\alpha_i}\\
&\times \prod_{i=1}^\ell\sum_{\nu_i=0}^{\lfloor(1+v)n\rfloor}
\Big(\frac{1}{(1+\nu_i)^{1/4}}+\frac{1}{(1+|\nu_i-u_i|)^{1/4}}\Big)
\Big(\frac{1}{(1+\nu_i)^{1/4}}+\frac{1}{(1+|\nu_i-v_i|)^{1/4}}\Big),
\end{align*}
where $u_i, v_i>0$ are some numbers.
Clearly, each of the last sums can be bounded by four sums of the form
\begin{align*}
\sum_{\nu_i=0}^{\lfloor(1+v)n\rfloor} \frac{1}{(1+|\nu_i-w_i|)^{1/4}(1+|\nu_i-z_i|)^{1/4}}
\le cn^{1/2}.
\end{align*}
This last estimate apparently holds independently of $w_i$ and $z_i$.
Estimate (\ref{est4-Lambda-n}) follows from above.

\smallskip

We are now in a position to complete the proof of (\ref{est-Lambda-n}).
Estimates (\ref{est2-Lambda-n}) and (\ref{est3-Lambda-n}) readily imply
\begin{equation}\label{est6-Lambda-n}
|\LL_n(x, y)|
\le \frac{cn^{|\a|+d}}{(1+n^{1/2}\|x-y\|)^{k+|\a|+d}},
\end{equation}
while by (\ref{est1-Lambda-n}) and (\ref{est4-Lambda-n})
we have for $1\le i \le \ell$
$$ 
|\LL_n(x, y)|
\le \frac{c n^{d/2}}{\prod_{i=1}^\ell (x_iy_i)^{\a_i+1/2}\prod_{i=
\ell+1}^d ( n^{-1})^{\a_i+1/2}
   (1+ n^{1/2}\|x-y\|)^{k-|\a|}}.
$$ 
Clearly, this estimate holds for an arbitrary permutation $i_1, i_2, \dots, i_d$
of the indices $1, 2, \dots, d$.
These estimates and (\ref{est6-Lambda-n}) yield
\begin{equation}\label{est7-Lambda-n}
|\LL_n(x, y)|
\le \frac{c n^{d/2}}{\prod_{i=1}^d(x_i y_i+n^{-1})^{\a_i+1/2}(1+n^{1/2}\|x-y\|)^{k-|\a|}}.
\end{equation}
To complete the proof 
we need the following simple inequality:
For $x,y\in \RR_+^d$
\begin{equation}\label{eq:px1}
(x_i+n^{-1/2})(y_i+n^{-1/2}) \le 3(x_iy_i+n^{-1})(1+n^{1/2}\|x-y\|), \quad 1\le i\le d.
\end{equation}
Combining these with (\ref{est7-Lambda-n}) we get
$$
|\LL_n(x, y)|
\le \frac{cn^{d/2}}{\prod_{i=1}^d(x_i+n^{-1/2})^{\a+1/2}(y_i+n^{-1/2})^{\a_i+1/2}
      (1+n^{1/2}\|x-y\|)^{k-2|\a|-d/2}},
$$
which implies (\ref{est-Lambda-n}) since $k$ was select so that
$k\ge \sigma+2|\a|+d/2$.

\smallskip

The proof of (\ref{est-Lambda-n-large}) is trivial.
Indeed, by Lemma~\ref{lem:asympt-Ln} it follows that
\begin{equation}\label{Est-F-large}
|\cF_n^\a(x)| \le cx^{-\a}e^{-\gamma x^2}
\quad\mbox{for \, $x\ge (6(1+v)n+3\alpha+3)^{1/2}$}.
\end{equation}
From this it easily follows that if $\max\{\|x\|^2, \|y\|^2\}\ge (6(1+v)n+3\alpha+3)^{1/2}$,
then
$$
|\Lambda_n(x, y)| \le cn^d e^{-\gamma \max\{\|x\|^2, \|y\|^2\}},
\quad \gamma>0,
$$
which readily implies (\ref{est-Lambda-n-large}).
$\qed$


\subsection{Proof of estimates (\ref{est-derivative}) and (\ref{est-derivative-large})
in Theorem \ref{thm:est-Lambda-n}}
Clearly, (\ref{est-derivative-large}) implies (\ref{est-derivative})
if $\max\{\|x\|, \|y\|\} \ge (6(1+v)n +3\a+3)^{1/2}$.

Assume $\max\{\|x\|, \|y\|\} < (6(1+v)n +3\a+3)^{1/2}\le cn^{1/2}$.
We will prove (\ref{est-derivative}) in this case by using the scheme of the proof of
(\ref{est-Lambda-n}) with appropriate modifications.
First, we need information about the derivative of $\PsiP_n^\a$.
The Laguerre polynomials satisfy the relation \cite[(5.1.14)]{Sz}
\begin{equation}\label{Laguerre-derivative}
    \frac{d}{dx} L_n^\a(x) = -L_{n-1}^{\a+1}(x) = x^{-1}
             \left[ n L_n^\a(x) - (n+\a) L_{n-1}^\a(x) \right].
\end{equation}
After taking the derivative of $\PsiP_n^\a$ (see (\ref{def.Psi-n})),
the first identity in \eqref{Laguerre-derivative} yields
\begin{equation} \label{Psi-diff2}
  \frac{d}{dx} \PsiP_n^\a(x) =  - x \left[ \PsiP_n^\a(x) +
       2 \sqrt{n} \PsiP_{n-1}^{\a+1}(x)\right],
\end{equation}
and from the second identity we similarly get
\begin{equation}\label{derivative-Psi}
x \frac{d}{dx}  \PsiP_n^\a(x) = -x^2 \PsiP_n^\a(x)
+ 2n \PsiP_n^\a(x) - 2 b_n \PsiP_{n-1}^\a,
\quad b_n:=\sqrt{n(n+\a)}.
\end{equation}
Here and in what follows we assume $\PsiP_k^\a(x) = 0$ for $k < 0$.
Also, from the recurrence relation for Laguerre polynomials \cite[(5.1.10)]{Sz}
one readily derives the identity
$$
x L_n^\a(x) = (2n +\a+1)L_n^\a(x) - (n+1) L_{n+1}^\a - (n+\a)L_{n-1}^\a(x),
\quad n\ge 1,
$$
with $L_0^\a (x) =1$ and $L_1^\a(x) = -x + \a +1$.
From this with the definition of
$\PsiP_n^\a$ in \eqref{def.Psi-n}, we get
\begin{equation} \label{Psi-recur}
  x^2 \PsiP_n^\a(x) = - b_{n+1}\PsiP_{n+1}^\a(x)
   + (2 n +\a + 1) \PsiP_n^\a(x) -  b_n \PsiP_{n-1}^\a(x),
\end{equation}
where $b_n$ is as above.
Combining this with (\ref{derivative-Psi}) gives
\begin{equation} \label{Psi-diff1}
   \frac{d}{dx} \PsiP_n^\a(x) = x^{-1} \left[-(\a +1)\PsiP_n^\a(x) + b_{n+1}
               \PsiP_{n+1}^\a(x)- b_n \PsiP_{n-1}^{\a}(x) \right].
\end{equation}
We also need the relation \cite[(5.1.13)]{Sz}
\begin{equation} \label{Laguerre-difference}
     L_n^\a(x) = L_{n}^{\a+1}(x) - L_{n-1}^{\a+1}(x).
\end{equation}
From this and (\ref{def.Psi-n}) we deduce
\begin{equation} \label{FFF}
 \PsiP_n^\a(x) = \sqrt{n+\a+1} \PsiP_n^{\a+1}(x) - \sqrt{n} \PsiP_{n-1}^{\a+1}(x).
\end{equation}
Using this identity with $\alpha$ replaced by $\alpha-1$,
\eqref{Psi-diff1}, and the obvious fact that $b_n = n+\CO(1)$,
we arrive at
\begin{equation} \label{Psi-diff-est}
  \left| \frac{d}{dx} \PsiP_n^\a(x)\right| \le c  x^{-1} \left[\max_{n-1\le m \le n+1}
          |\PsiP_m^\a(x)| + n^{1/2} \max_{n \le m \le n+1} |\PsiP_{m}^{\a-1}(x)| \right].
\end{equation}

By (\ref{Psi-diff2}) and (\ref{Est-norm-F})     
we readily get the estimate
$
\Big|\frac{d}{dx}\cF_n^\a(x)\Big|\le cxn^{\a/2+1}
$
and by (\ref{Psi-diff-est}) and (\ref{Est-norm-F})
$
\Big|\frac{d}{dx}\cF_n^\a(x)\Big|\le cx^{-1}n^{\a/2}.
$
Therefore,
\begin{equation}\label{Deriv-F}
\Big|\frac{d}{dx}\cF_n^\a(x)\Big|\le cn^{\a/2}\min\{x^{-1}, nx\}
\le cn^{(\a+1)/2},
\quad x\in\RR_+.
\end{equation}
We use this estimate to obtain
\begin{align}\label{partial-rep2}
\Big|\frac{\partial}{\partial x_r } \LL_n(x,y)\Big|
&\le \sum_{m=0}^{\lfloor(1+v)n\rfloor} \sum_{|\nu| =m}
|\PsiP_\nu^\a(y)|\Big|\frac{\partial}{\partial x_r }\PsiP_\nu^\a(x)\Big|\\
&\le c n^{1/2}\sum_{m=0}^{\lfloor(1+v)n\rfloor} \sum_{|\nu|=m} \nu^{\a}\le cn^{|\a|+d+1/2}.\notag
\end{align}

\smallskip

We next prove an analogue of (\ref{est4-Lambda-n}).
Let $0<x\le cn^{1/2}$.
Assuming that
$m\in \RR\setminus(x^2/5, x^2/3)$
we derive as before from
(\ref{Est-Fna1}) and (\ref{Psi-diff2})
\begin{align}\label{Est-dF1}
\Big|\frac{d}{dx}\cF_m^\a(x)\Big|
&\le x(|\cF_m^\a(x)|+2m^{1/2}|\cF_{m-1}^{\a+1}(x)|)\\
&\le cx\Big(\frac{1}{x^{\a+1/2}m^{1/4}}+\frac{m^{1/2}}{x^{\a+3/2}m^{1/4}}\Big)
\le \frac{cn^{1/2}}{x^{\a+1/2}m^{1/4}}. \notag
\end{align}
From (\ref{Est-Fna2}) and (\ref{Psi-diff2}) we similarly obtain
\begin{align}\label{Est-dF2}
\Big|\frac{d}{dx}\cF_m^\a(x)\Big|
\le \frac{cn^{1/2}}{x^{\a+1/2}(1+|4m-x^2|)^{1/4}}
\quad \mbox{for $m\in (x^2/5, x^2/3)$.}
\end{align}
We further proceed exactly as in the proof of (\ref{est4-Lambda-n}),
with $\cF_{\nu_r}^{\a_r}(x_r)$ replaced by
$\frac{\partial}{\partial x_r}\cF_{\nu_r}^{\a_r}(x_r)$
and for this term estimates (\ref{est-Fa-1})-(\ref{est-Fa-2}) are
replaced by (\ref{Est-dF1})-(\ref{Est-dF2}),
and we also use (\ref{Deriv-F}).
As a result, we get
\begin{equation}\label{partial-est-4}
\Big|\frac{\partial}{\partial x_r}\LL_n(x, y)\Big| \le \frac{cn^{(d+1)/2} }
{ \prod_{i=1}^\ell (x_iy_i)^{\a_i+1/2} \prod_{i=\ell+1}^d (n^{-1})^{\a_i+1/2}},
\quad 1\le \ell\le d.
\end{equation}

\smallskip


We now derive our main bound on $|(\partial/\partial x_r)\Lambda_n(x, y)|$.
It will be convenient to use the notation
$\partial f(t) : = f'(t)$.
After differentiating the expression of
$\LL_n(x,y)$ in \eqref{rep2-Ln} we obtain for $1 \le r \le d$,
\begin{equation} \label{Gn1+Gn2}
\frac{\partial}{\partial t_r} \LL_n(x,y) = \QQ_1(x,y) + \QQ_2(x,y),
\end{equation}
where
\begin{align}
& \QQ_1(x,y) :=
      \int_{[0,\pi]^d} \partial \KK_n^{k+|\a|+d}
        \Big(\|x\|_2^2 + \|y\|_2^2 +2 \sum_{i=1}^d x_i y_i \cos \theta_i \Big) \label{Gn1} \\
&  \qquad \times (2 x_r - 2y_r \cos \theta_r)
      \prod_{i=1}^d j_{\a_i-\frac12}(x_iy_i \cos\theta_i) (\sin\theta_i)^{2\a_i} d\theta,
       \notag\\
& \QQ_2(x,y) :  =
      \int_{[0,\pi]^d} \KK_n^{k+|\a|+d}
          \Big(\|x\|_2^2 + \|y\|_2^2 +2 \sum_{i=1}^d x_i y_i \cos \theta_i \Big) \label{Gn2}\\
 & \qquad \times  \prod_{i=1, \,i\ne r}^d j_{\a_i-\frac12}(x_iy_i \cos\theta_i)
                \partial j_{\a_r-\frac12}(x_r y_r \cos\theta_r)  y_r \cos\theta_r
                  (\sin\theta_i)^{2\a_i} d\theta . \notag
\end{align}
We first estimate $\QQ_1(x, y)$. By the left-hand-side identity in
\eqref{Laguerre-derivative} and \eqref{Laguerre-difference}
$$
\frac{d}{dt} [L_n^\a(t) e^{-t/2}]
= -(1/2)(L_n^\a(t) + 2L_{n-1}^{\a+1}(t))e^{-t/2}
= -(1/2)(L_n^{\a+1}(t)+ L_{n-1}^{\a+1}(t))e^{-t/2}.
$$
Hence, by the definition of $\KK_n^\lambda$ in (\ref{Def-KK}),
\begin{equation*} 
\partial \KK_n^{k+|\a|+d}(t)
= -[\KK_n^{k+|\a|+d+1}(t)+ \widetilde{\KK}_n^{k+|\a|+d+1}(t)]/2,
\end{equation*}
where $\widetilde{\KK}_n^\lambda(t)$ is define as $\KK_n^\lambda(t)$
but with $L_m^\lambda$ in the sum in (\ref{Def-KK}) replaced by $L_{m-1}^\lambda$.
Evidently, $\widetilde{\KK}_n^\lambda(t)$ has the same properties as
$\KK_n^\lambda(t)$.
Substituting the above in (\ref{Gn1}) and taking into account
(\ref{est-j-alpha})-(\ref{est-Fnk}) we get
$$
 \QQ_1(x,y) \le c n^{(-k + |\a| + d +1)/2} \int_{[0,\pi]^d}
    \frac{ |x_r - y_r \cos t_r|
      \prod_{i=1}^d t_i^{2\a_i} dt }
        {(\|x  - y\|_2^2 + \sum_{i=1}^d x_i y_i  t_i^2)^{(k+|\a| + d+1)/2} }.
$$
Now, using the fact that
$$
|x_r - y_r\cos t_r|
\le |x_r -y_r|+ 2x_r y_r \sin^2 (t_r/2)
\le |x_r -y_r|+ x_r^{-1}(x_r y_rt_r^2)
$$
and noticing that $|x_r -y_r|$ can be cancelled
by an $1/2$ power of the main term in the denominator, whereas $x_r y_rt_r^2$
needs a square of that much, we conclude that
\begin{align*}
   \QQ_1(x,y) & \le  c n^{(-k + |\a| + d +1)/2} \int_{[0,\pi]^d}
       \frac{ \prod_{i=1}^d t_i^{2\a_i} dt }
        {(\|x  - y\|_2^2 + \sum_{i=1}^d x_i y_i  t_i^2)^{(k+|\a| + d)/2}} \\
       & + cx_r^{-1}n^{(-k + |\a| + d +1)/2} \int_{[0,\pi]^d}
       \frac{ \prod_{i=1}^d t_i^{2\a_i} dt }
        {(\|x  - y\|_2^2 + \sum_{i=1}^d x_i y_i  t_i^2)^{(k+|\a| + d-1)/2} }.
\end{align*}
Both of the above integrals are of the form of $M_n^{k,\a}$ defined
in (\ref{eq:Jn}). In fact, we have
\begin{equation} \label{Gn1-est1}
\QQ_1(x,y) \le  c n^{1/2} M_n^{k,\a}(x,y) + cx_r^{-1}M_n^{k-1,\a}(x,y).
\end{equation}
Furthermore, evidently
$
|x_r - y_r \cos t_r| \le |x_r - y_r| + x_r t_r^2
$
and inserting $t_r^2$ into
the weight function of the integral, we obtain as above
\begin{align}\label{Gn1-est2}
\QQ_1(x,y)
&\le  c n^{1/2} M_n^{k,\a}(x,y) + cx_rM_n^{k,\a+e_r}(x,y).
\end{align}


We next estimate $\QQ_2$.
Using the integral representation \eqref{Bessel} for $j_{\a - \frac12}(x)$
we get
$$
\partial j_{\a - \frac{1}{2}} (x)  = c  \int_{-1}^1 e^{i xt} t (1-t^2)^{\a-1} dt,
\quad \a >0,
$$
while $\partial j_{-\frac12}(x) = c \sin x$.
Therefore,
$| \partial j_{\a-\frac12} (x) | \le c$ for $\a \ge 0$.
Consequently, using also that $y_r\le cn^{1/2}$, we obtain as in \eqref{eq:Jn}
\begin{align}\label{est-Gn2}
|\QQ_2(x,y)|
&\le cn^{1/2}M_n^{k,\a}(x,y).
\end{align}

Combining \eqref{Gn1-est1} and \eqref{est-Gn2} gives
\begin{equation}\label{partial-estA}
\left| \frac{\partial}{\partial x_r } \LL_n(x,y) \right|
\le cx_r^{-1} M_n^{k-1,\a}(x,y) +  c n^{1/2}M_n^{k,\a}(x,y),
\end{equation}
whereas combining \eqref{Gn1-est2} and \eqref{est-Gn2} gives
\begin{align} \label{partial-estB}
\left| \frac{\partial}{\partial x_r } \LL_n(x,y) \right|
\le cn^{1/2}M_n^{k,\a}(x,y) + c x_r M_n^{k, \a+e_r}(x,y).
\end{align}

We are now in a position to establish estimate \eqref{est-derivative}.
Using \eqref{est2-Lambda-n} in \eqref{partial-estA} and combining the
result with (\ref{partial-rep2}), we conclude that for $x_r \ge n^{-1/2}$
\begin{align}\label{partial-est-1}
 \left| \frac{\partial}{\partial x_r } \LL_n(x,y) \right|
        \le \frac{ c n^{|\a|+d+1/2}} {(1+n^{1/2}\|x-y\|)^{k+|\a|+d-1}}.
\end{align}
On the other hand, using \eqref{est2-Lambda-n} in \eqref{partial-estB} and
combining the result with (\ref{partial-rep2}) shows that estimate
(\ref{partial-est-1}) holds for $x_r \le n^{-1/2}$ as well.
Therefore, (\ref{partial-est-1}) holds for all $x, y \in \RR^d_+$.

In going further, using (\ref{est1-Lambda-n}) in \eqref{partial-estA} and combining
the result with \eqref{partial-est-4},
we obtain for $x_r \ge n^{-1/2}$ and $1 \le i \le \ell$
\begin{align} \label{partial-est-2}
 \left| \frac{\partial}{\partial x_r } \LL_n(x,y) \right|
\le \frac{c n^{(d+1)/2} } {\prod_{i=1}^\ell (x_iy_i)^{\a_i+\frac12}
         \prod_{i=\ell+1}^d (n^{-1})^{\a_i+\frac12}(1+ n^{1/2}\|x-y\|)^{k-|\a|-1}}.
\end{align}
On the other hand, using (\ref{est1-Lambda-n}) in \eqref{partial-estB}
and combining the result with \eqref{partial-est-4}, we see that the same
bound \eqref{partial-est-2} holds for $x_r \le n^{-1/2}$ as well.
Therefore, (\ref{partial-est-2}) holds in general.
Moreover, (\ref{partial-est-2}) holds for all possible permutations of the indices
and combining it with (\ref{partial-est-1}) leads to
\begin{align*} 
 \left| \frac{\partial}{\partial x_r } \LL_n(x,y) \right|
 \le \frac{c n^{(d+1)/2}} {\prod_{i=1}^d
     (x_iy_i+n^{-1})^{\a_i+\frac12} ( 1+ n^{1/2}\|x-y\|)^{k-|\a|-1}}.
\end{align*}
Now, estimate (\ref{est-derivative}) follows using (\ref{eq:px1}) as before.

\smallskip


The proof of (\ref{est-derivative-large}) is simple.
By (\ref{Est-F-large}) and (\ref{Psi-diff2}) it follows that
$$
\Big|\frac{d}{dx}\cF_n^\a(x)\Big|
\le cx^{-\a+1}e^{-\gamma x^2}
\le ce^{-\gamma' x^2}
\quad\mbox{for \, $x\ge (6(1+v)n+3\alpha+3)^{1/2}$}.
$$
This and (\ref{Est-F-large}) imply that if $\max\{\|x\|^2, \|y\|^2\}\ge (6(1+v)n+3\alpha+3)^{1/2}$,
then
$$
\left| \frac{\partial}{\partial x_r } \LL_n(x,y) \right|
\le cn^d e^{-\gamma'' \max\{\|x\|^2, \|y\|^2\}},
\quad \gamma''>0,
$$
which yields (\ref{est-derivative-large}).
$\qed$

\subsection{Proof of other localization estimates}\label{Other-loc-est}
${ }$

\smallskip

\noindent
{\bf Proof of Lemma~\ref{lem:instrument}.}
We will derive estimate (\ref{instrument}) from the following estimate:
If $s\in \RR$, $\gamma \ge 0$, $\sigma> (2\gamma+1)(|s|+1)+1$, and $z>0$, then
\begin{equation} \label{instrument1}
I:=\int_0^\infty \frac{u^{2\gamma+1}\, du}{(1+u)^{(2\gamma+1)s}(1+|u-z|)^\sigma}
\le \frac{c}{(1+z)^{(2\gamma+1)(s-1)}}.
\end{equation}

Consider first the case when $s\ge 1$.
Then
$I=\int_0^{z/2}+ \int_{z/2}^\infty=:J_1+J_2$.
Evidently,
$$J_1\le (1+z)^{-\sigma}\int_0^z1 du \le c(1+z)^{-\sigma+1}$$ and
$$
J_2\le \frac{c}{(1+z)^{(2\gamma+1)(s-1)}}\int_{z/2}^\infty\frac{du}{(1+|u-y|)^\sigma}
\le \frac{c}{(1+z)^{(2\gamma+1)(s-1)}} \qquad (\sigma>1).
$$
Since $\sigma > (2\gamma+1)(s-1)+1$ the above estimates for $J_1$ and $J_2$
yield (\ref{instrument1}).

Let $s<1$. Then we have
\begin{align*}
I&\le \int_0^\infty \frac{(1+u)^{(2\gamma+1)(1-s)}\, du}{(1+|u-z|)^\sigma}
= \int_{-z}^\infty \frac{(1+v+z)^{(2\gamma+1)(1-s)}\, du}{(1+|v|)^\sigma}\\
&\le c\int_{-\infty}^\infty \frac{(1+|v|)^{(2\gamma+1)(1-s)} + z^{(2\gamma+1)(1-s)}}
{(1+|v|)^\sigma}\, du\\
&\le c \int_{-\infty}^\infty \frac{du}{(1+|v|)^{\sigma+(2\gamma+1)(s-1)}}
+ cz^{(2\gamma+1)(1-s)}\int_{-\infty}^\infty \frac{du}{(1+|v|)^\sigma}\\
&\le \frac{c}{(1+z)^{(2\gamma+1)(1-s)}}.
\end{align*}
Here we used that $\sigma > (2\gamma+1)(1-s)+1$.
Therefore, (\ref{instrument1}) holds when $s<1$ as well.

We now proceed with the proof of (\ref{instrument}).
Denote by $J$ the integral in (\ref{instrument}).
Using that $|x_j-y_j|\le \|x-y\|$, we get
\begin{align*}
J &
\le \prod_{i=1}^d\int_0^\infty \frac{y_i^{2\a_i+1}dy_i}{(y_i+n^{-1/2})^{(2\a_i+1)s}
   (1+ n^{1/2}|x_i-y_i|)^{\sigma/d}}\\
&=n^{(2|\a|+d)s}\prod_{i=1}^d\int_0^\infty \frac{y_i^{2\a_i+1}dy_i}{(1+n^{1/2}y_i)^{(2\a_i+1)s}
   (1+ |n^{1/2}x_i-n^{1/2}y_i|)^{\sigma/d}}\\
&=n^{(2|\a|+d)(s-1)-d/2}\prod_{i=1}^d\int_0^\infty \frac{u^{2\a_i+1}du}{(1+u)^{(2\a_i+1)s}
   (1+ |u-n^{1/2}x_i|)^{\sigma/d}}\\
&\le cn^{(2|\a|+d)(s-1)-d/2}\prod_{i=1}^d\frac{1}{(1+n^{1/2}x_i)^{(2\a_i+1)(s-1)}}
=\frac{cn^{-d/2}}{\WW(n;x)^{s-1}}.
\end{align*}
Here for the last inequality we used (\ref{instrument1}).
$\qed$

\smallskip


\medskip\noindent
{\bf Proof of Theorems \ref{thm:est-tilde-Lambda-n} and \ref{thm:est-Lambda*-n}.}
By (\ref{rel-Psi-L}) we have
$\CL_\nu^\a(x)=2^{-1/2}\PsiP_\nu^\a(x^{\one/2})x^{\a/2}$
and by and (\ref{def.M-n})
$\cM^\a_\nu(x)=x^{\alpha+\one/2}\cF^\a_\nu(x)$
and hence
$$ 
\tLL_n(x, y)=2^{-1}\LL_n(x^{\1/2}, y^{\1/2}) x^{\a/2}y^{\a/2}
\quad\mbox{and}\quad
\LL_n^*(x, y)=\LL_n(x, y) x^{\a+\one/2}y^{\a+\one/2}
$$ 
Now, it is easy to see that these relations and estimates (\ref{est-Lambda-n})
and (\ref{est-derivative}}) yeild
(\ref{est-tilde-Lambda-n}) and (\ref{est-tilde-derivative}) as well as
(\ref{est-Lambda*-n}) and (\ref{est-derivative*}).
$\qed$


\subsection{Proof of Lemma~\ref{lem:lower-bound}}

The main step is to prove Lemma \ref{lem:lower-bound} for dimension $d=1$.
To this end we will need a lemma which goes back to van der Corput
(see e.g. \cite[Vol. I, p. 197-198]{Z}).

\begin{lem}\label{Corput}
If $f''(u) \ge \rho > 0$ or $f''(u) \le - \rho <0$ on $[a, b]$, then
$$
\Big| \sum_{a \le n \le b} e^{2 \pi i f(n)} \Big|
\le (|f'(b) - f'(a)| + 2)
(4 \rho^{-1/2} + c).
$$
\end{lem}

Evidently, when $d=1$ Lemma \ref{lem:lower-bound} is immediate from the following lemma.


\begin{lem} \label{lem:low-bound-d=1}
For any $\eps >0$ and $\delta>0$ there exists a constant $c > 0$ such that for $n\ge 1/\eps$
\begin{equation} \label{lowerbd-d=1}
\AA_n(x) :=
e^{-x}\sum_{m=n}^{n + \lfloor\eps n\rfloor}  \frac{[L_m^\a(x) ]^2}{L_m^\a(0)}
  \ge c  n^{1/2} (x + \tfrac{1}{n})^{-\a-1/2}, \quad 0\le x \le (4-\delta)n.
\end{equation}
\end{lem}

\noindent
{\bf Proof.}
We may assume that $\eps\le 1$ and $n\ge n_0$, where $n_0$ is sufficiently large.
The proof uses the asymptotic of $L_n^\a(x)$ and is divided
into several cases.

\medskip

{\it Case 1:} $0 \le x < c^\diamond n^{-1}$ with $c^\diamond:= (\a+1)(\a+3)$
($c^\diamond n^{-1}$ is larger than the smallest zero of $L_n^\a$
\cite[(6.31.12)]{Sz}).
We need the asymptotic formula \cite[(8.22.4)-(8.22.5)]{Sz}
$$
e^{-x/2} x^{\alpha/2} L_n^\a(x)= N^{-\alpha} \frac{\Gamma(n+\alpha+1)}{n!}
J_\alpha(2 (Nx)^{1/2}) + x^{\a/2+2}\CO(n^\alpha),
\quad 0<x\le c/n,
$$
where $N = n+(\alpha+1)/2$.
Using also that
$
J_\alpha(z) =\frac{z^\a}{2^\a\Gamma(\a+1)}+\CO(z^{\a+2}),
$
we obtain
$$
e^{-x/2} L_n^\a(x) \sim n^{\alpha} +  x^2\CO(n^{\alpha})
\ge c n^{\alpha},
\quad 0 \le x < c/n.
$$
Combining this with
$L_n^\alpha(0) = \binom{n+\alpha}{n} \sim n^\alpha$
we arrive at
$$
\AA_n(x) \ge c \sum_{m=n}^{n+\lfloor \varepsilon n \rfloor} m^\a \sim n^{\alpha+1},
\quad 0 \le x < c^\diamond n^{-1}, 
$$
which proves (\ref{lowerbd-d=1}) in this case.

\medskip

{\it Case 2:} $c^\diamond n^{-1} \le x \le c_* n^{-1}$, where the constant $c_* > 1$
will be selected later on.
In this case we use relation \eqref{lowerOP} and \eqref{zeros-difference} to conclude
that
$$
    e^{-x} |L_n^\a(x) |^2 \sim n^{2\a+2} (x - t_{k_x,n})^2.
$$
Furthermore, by a theorem of Tricomi (see \cite{Ga} for the references),
we know that for all the zeros of $L_n^\a$ in the interval $0 < x < c/n$
we have
$
t_{k,n} = \frac{j_{\a,k}^2}{n} (1+ \CO(n^{-2}))
$
as $n \to \infty$, where $j_{\a,k}$, $k = 1,2,...$, are the positive zeros, in
increasing order, of the Bessel function $J_\a(x)$. Consequently,
\begin{align*}
   \AA_n(x) &  \ge c n^\a \sum_{m=n}^{n+\lfloor\eps n\rfloor}
      \left( (m x - j_{\a, k_x}^2)^2 - c m^{-1} |m x  - j_{\a,k_x}|\right) \\
         &  \ge c n^\a \biggl( \sum_{m=n}^{n+\lfloor\eps n\rfloor}
             (m x - j_{\a, k_x}^2)^2 - c \biggr) \ge c n^{\a+1}.
\end{align*}
Here for the last estimate we used that $j_{\a,k}\to \infty$ as $k\to\infty$
and hence there are only finitely many zeros of $J_\a(x)$ such that
$j_{\a,k}^2\le c_* n^{-1}(n+\lfloor\eps n\rfloor) \le c$;
the argument is the same as in the analogous
situation for Jacobi polynomials in \cite{KPX1}.

\medskip

{\it Case 3:} $c_* n^{-1} \le x \le c^*$, where $c_*$ is sufficiently large and its
value is to be determined.
In this case we use the asymptotic formula for $L_n^\a(x)$ \cite[(8.22.6)]{Sz}:
\begin{align*}
e^{-x/2} L_n^\a(x) = & \pi^{-1/2} x^{-\a/2-1/4} n^{\a/2-1/4} \\
& \times \left[ \cos \left( 2(nx)^{1/2} - \a\pi/2 - \pi/4 \right )
          + \CO(1)(n x )^{-1/2} \right],
\end{align*}
which holds for $c'n^{-1}\le x \le c''$ and $\CO(1)$ depends only on $c', c''$.
We denote $\g := \a \pi/2 + \pi/4$ and deduce from above
\begin{align*}
x^{\a+1/2} \AA_n(x)
&\ge c n^{-\a} e^{-x} x^{\a+1/2}
     \sum_{m=n}^{n + \lfloor\eps n\rfloor} \left[L_m^\a(x)\right]^2\\
&\ge c\sum_{m=n}^{n + \lfloor\eps n\rfloor} m^{-1/2}
           \left(\cos[2(mx)^{1/2} - \g]  + \CO(1)(n x )^{-1/2} \right)^2\\
&\ge c n^{-1/2} \sum_{m=n}^{n + \lfloor\eps n\rfloor}\cos^2[2(mx)^{1/2} - \g]
        + \CO(1) c_*^{-1/2} n^{1/2}.
\end{align*}
Using the fact that $2 \cos^2 t  = 1 + \cos 2 t$ and then write $2 \cos 2t =
  e^{2i t} + e^{-2 it}$, we see that
\begin{align*}
\Sigma
:= 4 \sum_{m=n}^{n + \lfloor\eps n\rfloor}\cos^2[2(mx)^{1/2} - \g]
\ge 2\lfloor\eps n\rfloor +  \sum_{m=n}^{n + \lfloor\eps n\rfloor}
\left [ e^{2\pi i(y\sqrt{m}-\g')} + e^{-2\pi i(y\sqrt{m} -\g')}\right],
\end{align*}
where $y :=(2/\pi)\sqrt{x}$ and $\g':=2\g/\pi$.
The last sum can be estimated by making use of Lemma \ref{Corput} with
$f (u) = y \sqrt{u}$, $a = n$ and $b= n + \lfloor\eps n\rfloor$.
We get
\begin{align*}
\Sigma
&\ge 2\lfloor\eps n\rfloor - 2(2 + x^{1/2}n^{-1/2})(c + 24x^{-1/4} n^{3/4})\\
&\ge 2\lfloor\eps n\rfloor - 2(2 + (c^*)^{1/2}n^{-1/2})(c + 24c_*^{-1/4} n).
\end{align*}
Putting the above estimates together, we arrive at
$$
x^{\a+1/2} \AA_n(x)
\ge  c n^{-1/2}
\left(2\lfloor\eps n\rfloor - 2(2 + (c^*)^{1/2}n^{-1/2})(c + 24c_*^{-1/4} n) \right)
       + \CO(1)  c_*^{-1/2} n^{1/2}.
$$
Choosing $c_*$ sufficiently large  shows that the right-hand side
of the above inequality is bounded below by $c n^{1/2}$ for sufficiently large $n$.
Thus (\ref{lowerbd-d=1}) is proved in this case.

\medskip

{\it Case 4:} $c^* \le x \le (4-\delta)n$. Here we apply another asymptotic formula
of Laguerre polynomial \cite[(8.22.9)]{Sz}:
For $x = (4m+2\a+2)\cos^2 \phi$
with $\eps \le \phi \le \pi/2 - \eps m^{-1/2}$,
\begin{align*}
x^{\a/2+1/4}e^{-x/2} L_m^\a(x)
& = (-1)^m (\pi \sin \phi)^{-1/2} m^{\a/2-1/4}\\
& \times \left\{\sin \Big[(m+ \tfrac{\a+1}{2})(\sin 2 \phi - 2\phi) + 3 \pi/4\Big] +
\CO(1) (m x )^{-1/2}\right\}.
\end{align*}
Note that the range of $x$ above covers the range of this case.
From above, as in Case 3, we obtain
\begin{align*}
& x^{\a+1/2} \AA_n(x)  \ge  c n^{-\a} e^{-x} x^{\a+1/2}
\sum_{m=n}^{n + \lfloor\eps n\rfloor} \left[L_m^\a(x)\right]^2\\
&  \qquad\qquad \ge  c n^{-1/2} \sum_{m=n}^{n + \lfloor\eps n\rfloor}
\sin^2[(m+\tfrac{\a+1}{2})(\sin 2 \phi - 2 \phi) + 3 \pi /4] + \CO(1)(c^*)^{-1/2}.
\end{align*}
The last sum is again bounded below by $c n$, which can be proved either
by using Lemma \ref{Corput} or by summing up using simple trigonometric identities.
This shows again that (\ref{lowerbd-d=1}) holds.
$\qed$

\medskip

\noindent
{\bf Proof of (\ref{lower-bound}) in the case \boldmath $d\ge 2$}.
We may again assume $\eps \le 1$.
We will use induction on $d$. To indicate the dependence of $\cF_m^\a$ on $d$
we write $\cF_{m,d}^\a: =\cF_m^\a$. Assume that \eqref{lowerbd} has been
established for dimensions up to $d-1$. By definition
$$
\cF_{m,d}^\a(x,x) = \sum_{k=0}^m \Big[\cF_k^{\a_d}(x_d)\Big]^2\cF_{m-k,\,d-1}^{\a'}(x',x'),
\quad x = (x',x_d),\;\: \alpha=(\alpha', \alpha_d),
$$
and hence
\begin{align}\label{inductive-step}
\sum_{m=n}^{n+\lfloor d\eps n\rfloor}\cF_{m,d}^\a(x,x)
&\ge \sum_{m=n}^{n+\lfloor d\eps n\rfloor}\sum_{k=0}^{\lfloor\eps n\rfloor}
\Big[\cF_k^{\a_d}(x_d)\Big]^2\cF_{m-k,\,d-1}^{\a'}(x',x')\notag\\
&=\sum_{k=0}^{\lfloor\eps n\rfloor}\Big[\cF_k^{\a_d}(x_d)\Big]^2
\sum_{m=n}^{n+\lfloor d\eps n\rfloor}\cF_{m-k,\,d-1}^{\a'}(x',x')\\
&\ge \sum_{k=0}^{\lfloor\eps n\rfloor}\Big[\cF_k^{\a_d}(x_d)\Big]^2
\sum_{j=n}^{n+\lfloor (d-1)\eps n\rfloor}\cF_{j,\,d-1}^{\a'}(x',x').\notag
\end{align}
It follows by (\ref{def.Psi-n}) and (\ref{def-Kn})-(\ref{Christoph}) that
for $0\le x\le \sqrt{(4-\delta)n}$
$$
\sum_{k=0}^n\Big[\cF_k^{\a}(x)\Big]^2=ce^{-x^2}K_n^\a(x^2, x^2)
\ge cn^{1/2}\Big(x^2+\frac 1n\Big)^{-\alpha-1/2}
\ge cn^{1/2}(x+n^{-1/2})^{-2\alpha-1}.
$$
Combining this estimate with (\ref{inductive-step}) and the inductive assumption shows that
(\ref{lowerbd}) holds in dimension $d$.
$\qed$

\bigskip


\noindent
{\bf Proof of Proposition \ref{prop:norms}.}
We first prove (\ref{norms2}).
Let $g\in V_n$.
Assume $1<q<\infty$ and let $\LL_n$ be the kernel from (\ref{def.Ln}),
with $\ha$ admissible of type $(a)$.
Evidently $g=\LL_n*g$ and using H\"older's inequality and
Proposition~\ref{prop:Lp-bound} we obtain for $x\in \RR^d_+$
\begin{align*}
&|g(x)|
\le \|\WW(n;\cdot)^{\ss +\frac{1}{p} - \frac{1}{q}} g(\cdot)\|_q
\left( \int_{\RR^d_+} \Big|\LL_n(x,y) \WW(n;y)^{-\ss-\frac{1}{p}+
\frac{1}{q}}\Big|^{q'} \w(y) dy \right)^{\frac{1}{q'}} \\
& \le c \frac{n^{d/2}} {\WW(n;x)^{1/2}}
\left (\int_{\RR^d_+}
\frac{\w(y)dy}{ \WW(n;y)^{\frac{q'}{2} + \beta} (1+ n^{1/2}\|x-y\|)^\sigma }\right)^{\frac{1}{q'}}
\|\WW(n;\cdot)^{\ss +\frac{1}{p} - \frac{1}{q}} g(\cdot)\|_q,
\end{align*}
where $\beta := q'(\ss + \frac{1}{p}-\frac{1}{q})$.
To estimate the last integral we use estimate (\ref{instrument}) from
Lemma~\ref{lem:instrument} to obtain
\begin{equation} \label{infinity-q}
 |g(x)|  \le  c \frac{n^{d/2q}}{\WW(n;x)^{\ss +1/p}}
     \|\WW(n;\cdot)^{\ss +\frac{1}{p} - \frac{1}{q}} g(\cdot)\|_q
\end{equation}
and hence
\begin{equation} \label{aaa2}
\|\WW(n;\cdot)^{\ss+\frac{1}{p}}g(\cdot)\|_\infty \le c n^{d/2q}
\|\WW(n;\cdot)^{\ss +\frac{1}{p} - \frac{1}{q}} g(\cdot)\|_q,
\quad 1<q\le \infty.
\end{equation}
If $0 < q \le 1$, then the above estimate with $q=2$ gives
\begin{align*}
 \|\WW(n;\cdot)^{\ss+\frac{1}{p}}g(\cdot)\|_\infty
&\le c n^{d/4}
 \|\WW(n;\cdot)^{\ss +\frac{1}{p} - \frac{1}{2}} g(\cdot)\|_2 \\
 & \le c n^{d/4}  \| \WW(n;\cdot)^{\ss+1/p}g(\cdot)\|_\infty^{1-q/2}
 \|\WW(n;\cdot)^{\ss +\frac{1}{p} - \frac{1}{q}} g(\cdot)\|_q^{q/2}.
\end{align*}
Consequently, (\ref{aaa2}) holds for $0 < q \le 1$ as well.

Let $0<q<p<\infty$.
Using (\ref{aaa2}), we have
\begin{align*}
&\|\WW(n;\cdot)^\ss g(\cdot)\|_p\\
&\qquad = \left( \int_{\RR^d_+} \left|\WW(n;x)^{\ss+\frac{1}{p}} g(x)\right|^{p-q}
\left|\WW(n;x)^{\ss+\frac{1}{p}-\frac{1}{q}} g(x)\right|^{q}\w(x) dx \right)^{1/p} \\
&\qquad \le \|\WW(n;\cdot)^{\ss+\frac{1}{p}}g(\cdot)\|_\infty^{1-q/p}
\|\WW(n;\cdot)^{\ss +\frac{1}{p} - \frac{1}{q}} g(\cdot)\|_q^{q/p}\\
&\qquad = c n^{(d/2)(1/q-1/p)}
\|\WW(n;\cdot)^{\ss +\frac{1}{p} - \frac{1}{q}} g(\cdot)\|_q.
\end{align*}
Hence (\ref{norms2}) holds when $p < \infty$.
In the case $p = \infty$ (\ref{norms2}) follows from (\ref{aaa2}).

\smallskip


To prove (\ref{norms1}) we first assume that $1<q<\infty$.
We use again that $g=\LL_n*g$, H\"older's inequality,
Proposition~\ref{prop:Lp-bound},
and that $\WW(n;x)\ge n^{-|\alpha|-d/2}$
to obtain
$$
|g(x)| \le c\|g\|_q\left(n^{d/2}\WW(n;x)^{-1}\right)^{1/q}
\le c n^{(d+|\alpha|)/q}\|g\|_q,
\quad x\in \RR^d_+,
$$
and hence
$
\|g\|_\infty \le c n^{(d+|\alpha|)/q} \|g\|_q.
$
For the rest of the proof of (\ref{norms1}) one proceeds similarly as in
the proof of (\ref{norms2}). We skip the details.

\smallskip

To prove estimate (\ref{norms3}) we first observe that
(\ref{infinity-q}) with  $s=\gamma + 1/p -1/q$ yields
$$
|g(x)|  \le c \frac{n^{d/2q}} {W_\a(n;x)^{s+ 1/q}} \|W_\a(n;\cdot)^s g(\cdot)\|_q,
\quad 1 < q < \infty,
$$
and, since $W_\a(n; x) \ge n^{-|\a|-\frac{d}2}$, we get
$
\|g\|_\infty  \le c n^{(|\a|+\frac d 2)s + (|\a|+d)/q} \|W_\a(n;\cdot)^s g(\cdot) \|_q.
$
The remaining part of the proof is similar to the proof of (\ref{norms2}). We omit it.
$\qed$


\medskip

\noindent
{\bf Proof of Lemma~\ref{lem:char-S}.}
(a) By (\ref{est2-eLn}) and the definition of $\cF_n^\a$, it follows that
$\|\cF_n^\a\|_\infty \le c n^{\alpha/2}$.
Hence, using (\ref{Psi-diff1}) if $|x| \le 1$ and (\ref{Psi-diff-est})
if $|x| \ge 1$, we obtain
$$
    \left | \frac{d}{dx}\cF_n^\a(x) \right | \le c  n^{( \a+1)/2},
            \quad x \in \RR_+.
$$
Furthermore, taking one more derivative of (\ref{Psi-diff2})
and using (\ref{Psi-diff1}) shows that
\begin{align*}
 \frac{d^2}{d x^2} \cF_n^\a(x) & = - [\cF_n^\a(x) + 2 \sqrt{n} \cF_{n-1}^{\a+1}(x)]
    + x \frac{d}{d x} \cF_n^\a(x)+ 2 \sqrt{n}\, x \frac{d}{d x} \cF_{n-1}^{\a+1}(x) \\
  & =  - [\cF_n^\a(x) + 2 \sqrt{n} \cF_{n-1}^{\a+1}(x)]\\
   &\quad     -(\a+1) \cF_n^\a(x) + b_{n+1}\cF_{n+1}^\a(x) -  b_n \cF_{n-1}^\a(x) \\
   &\quad   + 2 \sqrt{n}\left [ -(\a+1) \cF_{n-1}^{\a+1}(x) + b_{n}\cF_{n}^{\a+1}(x) -
         b_{n-1} \cF_{n-2}^{\a+1}(x) \right],
\end{align*}
which allows us to iterate and express  $\frac{d^{k+1}}{d x^{k+1}} \cF_n^\a(x)$ in
terms of $\frac{d^{k-1}}{d x^{k-1}} \cF_n^\a(x)$  and $\frac{d^{k-1}}{d x^{k-1}}
\cF_n^{\a+1}(x)$. The recurrence relation (\ref{FFF})
allows us to use induction to conclude that
$$
  \left | \frac{d^k}{dx^k}\cF_n^\a(x) \right | \le c  n^{(\a+k)/2},
            \quad x \in \RR_+.
$$
Therefore, for the product Laguerre functions, we have
$$
  \left |\DD^\b \cF_\nu^\a(x)\right| \le c  (|\nu|+1)^{(|\a|+|\b|)/2}, \quad |\nu| =n,
          \quad \b \in \NN_0^d,  \quad x \in \RR_+^d.
$$
Furthermore, together with the three term relation (\ref{Psi-recur}),
the above inequality also shows that
$$
\left |x^{2 \g} \DD^\b \cF_\nu^\a(x)\right| \le c  (|\nu|+1)^{(|\a|+|\b|+2|\g|)/2},
    \quad |\nu| =n, \quad \b,\g \in \NN_0^d, \quad x \in \RR_+^d.
$$
Hence, if $|\langle \phi, \cF_\nu^\a \rangle | \le c_k(|\nu| +1)^{-k}$ for all $k$,
then
$$
x^{\g} \DD^\b \phi(x) = \sum_{\nu \in \NN_0^d}\langle \phi, \cF_\nu^\a \rangle
     x^{\g} \DD^\b \cF_\nu^\a(x),
$$
where the series converges uniformly and hence
\begin{equation}\label{Pab<P*}
 \left| x^{\g} \DD^\b \phi(x)\right| \le c \sum_{\nu \in \NN_0^d} \left|
     \langle \phi, \cF_\nu^\a \rangle \right|  (|\nu|+1)^{(|\a|+|\b|+2|\g|)/2}
       \le c_kP_k^*(\phi)
\end{equation}
if $k > d+ |\a|+|\b|+2|\g|)/2$, which shows that $\phi \in \CS_+$.

\smallskip


(b)
Assuming that $\phi \in \CS_+$ we next show that
$|\langle \phi, \cF_\nu^\a \rangle|$ has the claimed decay.
From the well
known second order differential equation satisfied by $L_n^\a$, a
straightforward computation shows that $\cF_n^\a(x)$ satisfies the equation
$$
    y''+ \frac{2\a+1}{x} y' - x^2 y + 2(2n+\a+1) y =0.
$$
In particular, it follow that $\cF_\nu^\a(x)$ satisfies, for each $i =1,2,\ldots,d$,
the equation
\begin{equation} \label{PDE}
\CD_{x_i} u   + x_i^2 u = 2( 2  \nu_i + \a_i +1)u, \quad
\hbox{where}\quad
\CD_{x_i} : =  -\partial_i^2 - (2 \a_i +1)x_i^{-1}\partial_i
\end{equation}
and
$\partial_i =\frac{\partial}{\partial x_i}$.

Let $k\ge 1$ and assume that the multi-index $\nu$ is fixed and
$\|\nu\|=\max_{1\le j \le d} \nu_j\ge k$.
Choose $i$ so that $\nu_i = \|\nu\|$
and denote
 $\wh x_i = (x_1,\ldots,x_{i-1}, 0, x_{i+1},\ldots,x_d)$.
Denote briefly
$
\UU_r(x) := \partial_i^r(\phi(x)e^{x_i^2/2}).
$
Then by Taylor's identity
$$
\phi(x)e^{x_i^2/2} -\sum_{r=0}^{2k-1}x_i^r\UU_r(\wh x)/r!
= \frac{x_i^{2k}}{(2k-1)!}\int_0^1(1-t)^{2k-1}\UU_{2k}(\wh x+tx_ie_i)dt,
$$
which easily leads to
\begin{align}\label{Taylor}
\phi_i(x)
&:= \phi(x) -e^{-x_i^2/2}\sum_{r=0}^{2k-1}x_i^r\UU_r(\wh x)/r!\\
&=x_i^{2k}\int_0^1(1-t)^{2k-1}\sum_{j=0}^{2k} b_{2k-j}(tx_i)
\partial_i^j\phi(\wh x + tx_ie_i)e^{-x_i^2(1-t^2)}dt,\notag
\end{align}
where $b_{j}(\cdot)$ $(0\le j\le 2k)$ is a polynomial of degree $\le j$
and $e_i$ is the $i$th coordinate vector in $\RR^d$.
Then by the orthogonality of $\cF_{\nu_i}^{\a_i}$ (recall that $\nu_i\ge 2k$)
and (\ref{PDE}) it follows that
\begin{align*}
\langle \phi, \cF_\nu^\a \rangle  =   \langle  \phi_i, \cF_\nu^\a \rangle
= \frac{1}{2( 2 \nu_i+ \a_i+1)}
\langle \phi_i, (\CD_{x_i}+x_i^2) \cF_\nu^\a \rangle.
\end{align*}
The operator $\CD_{x_i}$ can be written in a self-adjoint
form
$
    x_i^{2\a_i+1} \CD_{x_i} = \partial_i \left( x_i^{2\a_i+1} \partial_i \right).
$
We use this and integration by parts to obtain
\begin{align*}
\langle \phi_i, \CD_{x_i} \cF_\nu^\a \rangle
& = \int_{\RR^{d-1}_+} \int_{\RR_+}  \phi_i(x) \partial_i \left(  x_i^{2\a_i+1} \partial_i
       \cF_\nu^\a(x)\right)  dx_i d \wh x\\
& = \int_{\RR^{d-1}_+} \int_{\RR_+}  \partial_i \left(  x_i^{2\a_i+1} \partial_i \phi_i(x) \right)
           \cF_\nu^\a(x) dx_i d \wh x
     =  \langle  \CD_{x_i}\phi_i,   \cF_\nu^\a \rangle.
\end{align*}
Consequently,
\begin{align}\label{IterateCoeff}
\langle \phi, \cF_\nu^\a \rangle
&= \frac{1}{2( 2 \nu_i+ \a_i+1)}\langle (\CD_{x_i}+x_i^2)\phi_i, \cF_\nu^\a \rangle\\
&= \frac{1}{2^k( 2 \nu_i+ \a_i+1)^k}\langle (\CD_{x_i}+x_i^2)^k\phi_i, \cF_\nu^\a \rangle,\notag
\end{align}
where we iterated $k$ times.
It is easy to see that there is a representation of the form
$$
(\CD_{x_i}+x_i^2)^k
= (-\partial_i^2-(2 \a_i +1)x_i^{-1}\partial_i+x_i^2)^k
=\sum_{j=0}^{2k}\sum_{\ell=-2k}^{2k-j}a_{j\ell}x_i^{-\ell}\partial_i^j
$$
for some constants $a_{j\ell}$.
On the other hand, by (\ref{Taylor}) it follows that
if $j+\ell\le 2k$
$$
\sup_x |x_i^{-\ell}\partial_i^j\phi_i(x)|
\le c\max_{|\gamma|\le 4k, |\beta|\le 2k+j}\sup_x |x^\gamma\partial^\beta \phi(x)|
=c\max_{|\gamma|\le 4k, |\beta|\le 2k+j}P_{\beta, \gamma}(\phi).
$$
We use the above in (\ref{IterateCoeff}) to obtain
\begin{align}\label{coef<norm}
|\langle \phi, \cF_\nu^\a \rangle|
&\le \frac{1}{2^k(2 \nu_i+ \a_i+1)^k}
\max_{|\gamma|\le 4k, |\beta|\le 4k}P_{\beta, \gamma}(\phi)\|\cF_\nu^\a\|_1\\
&\le c|\nu|^{-k+(|\a|+d)/2}\max_{|\gamma|\le 4k, |\beta|\le 4k}P_{\beta, \gamma}(\phi),
\qquad \|\nu\|\ge k.\notag
\end{align}
Here we also used that $\|\cF_\nu^\a\|_1\le c|\nu|^{(|\a|+d)/2}$
which follows from Lemma~\ref{lem:asympt-Ln}.
Estimate (\ref{coef<norm}) shows that
$|\langle \phi, \cF_\nu^\a \rangle | \le c_k (|\nu|+1)^{-k+(|\a|+d)/2}$
for any $k\ge 1$.
Thus $|\langle \phi, \cF_\nu^\a \rangle |$ has the claimed decay.

The equivalence of the topologies on $\cS_+$ induced by the semi-norms
$P_{\gamma, \beta}$ from (\ref{schwartz}) and the norms $P_k^*$ from (\ref{def-P*})
follows readily by (\ref{Pab<P*}) and (\ref{coef<norm}).
$\qed$

\section{Proofs for Sections \ref{Tri-Liz}-\ref{Besov}\label{Proofs2}}
\setcounter{equation}{0}


\noindent
{\bf Proof of Proposition~\ref{prop:multipliers}.}
We shall use a standard decomposition of unity argument.
Suppose $\hb\in C^\infty(\RR)$ satisfies the conditions: $\supp \hb \subset [1/4, 4]$,
$b\ge 0$, and $\hb(t)+\hb(4t)=1$ on $[1/4, 1]$; hence
$\sum_{\ell=0}^\infty \hb(4^{-\ell}t)=1$, $t\in [1, \infty)$.
Now, define
$$
\Phi_0(x, y):= m(0)\cF_0^\a(x, y)
\quad\mbox{and}\quad
\Phi_\ell(x, y):=\sum_{j=0}^{4^\ell} m(j)\hb(j/4^{\ell-1})\cF_j^\a(x, y), \;\; \ell\ge 1.
$$
Then for the kernel $K(x, y)$ of the operator $T_m^\a$ we have
$
K(x, y)=\sum_{\ell=0}^\infty \Phi_\ell(x, y).
$
By (\ref{cond-on-m}) it readily follows that
$\|(d/dt)^k[m(t)\hb(t/4^{\ell-1})]\|_\infty \le c4^{-\ell k}$
and just as in the proof of Theorem~\ref{thm:est-Lambda-n} (using also (\ref{W<W}))
we get for $x, y\in\RR^d_+$
$$
|\Phi_\ell(x, y)|
\le \frac{c2^{\ell d}}{\W(4^\ell;y)(1+2^\ell\|x-y\|)^\sigma},
\;\;
\Big|\frac{\partial}{\partial y_r}\Phi_\ell(x, y)\Big|
\le \frac{c2^{\ell(d+1)}}{\W(4^\ell;y)(1+2^\ell\|x-y\|)^\sigma},
$$
for $1\le r\le d$, where $\sigma = k-(5/2)|\a|-(3/4)d-2$.
By a simple standard argument these two estimates ($\sigma >d+1$)
lead to
$$
|K(x, y)|\le \frac{c}{\w(y)\|x-y\|^d}
\quad\mbox{and}\quad
\Big|\frac{\partial}{\partial y_r}K(x, y)\Big|\le \frac{c}{\w(y)\|x-y\|^{d+1}},
\quad 1\le r\le d.
$$
As in the weighted case on $\RR^d$ (see \cite{Stein}),
these estimates show that $T_m^\a$ is a Calder\'{o}n-Zygmund type operator
and hence $T_m^\a$ is bounded on $L^p(\w)$, $1<p<\infty$.
$\qed$

\medskip


\noindent
{\bf Proof of Lemma $\ref{lem:Phi*psi}$.}
Using the orthogonality of Laguerre functions, we have
$\Phi_j\ast \psi_\xi(x)=0$ for $\xi\in\cX_m$ if $|m-j|\ge 2$.

Let $\xi \in \cX_m$, $j-1\le m \le j+1$.
Assume first that $\|\xi\|\le (1+\delta)\sqrt{6}\cdot2^m$.
From (\ref{local-Needlets1})-(\ref{local-needlets2}) it follows that
\begin{align*}
|\Phi_j\ast \psi_\xi(x)|
&\le  c_\sigma \frac{2^{m3d/2}}{\sqrt{\WW(4^m;x)}}\int_{\RR^d_+}\frac{\w(y)}
{\WW(4^m;y)(1+2^m\|x-y\|)^\sigma (1+2^m\|y-\xi\|)^\sigma}\, dy\\
&\le \frac{c2^{m3d/2}}{\sqrt{\WW(4^m;x)}}\int_{\RR^d}\frac{dy}
{(1+2^m\|x-y\|)^\sigma (1+2^m\|y-\xi\|)^\sigma}\qquad (\sigma>d)\\
&\le\frac{c2^{md/2}}{\sqrt{\WW(4^m;x)}(1+2^m\|x-\xi\|)^\sigma}
\le\frac{c2^{md/2}}{\sqrt{\WW(4^m;\xi)}(1+2^m\|x-\xi\|)^{\sigma-2|\alpha|-2d}}\\
&\le\frac{c}{\mu(R_\xi)^{1/2}(1+2^m\|x-\xi\|)^{\sigma-2|\alpha|-2d}},
\end{align*}
where for the last two inequalities we used (\ref{size-R-xi1})-(\ref{W<W}).
Since $\sigma$ can be arbitrarily large the claimed estimate (\ref{Phi*psi})
follows.

Let $\|\xi\|>(1+\delta)\sqrt{6}\cdot2^m$.
Just as above we use (\ref{local-Needlets1}) and (\ref{local-needlets3}) to obtain
\begin{align*}
|\Phi_j\ast \psi_\xi(x)|
&\le  c_\sigma \frac{2^{m(d-L)}}{\sqrt{\WW(4^m;x)}}\int_{\RR^d_+}\frac{\w(y)}
{\WW(4^m;y)(1+2^m\|x-y\|)^\sigma (1+2^m\|y-\xi\|)^\sigma}\, dy\\
&\le\frac{c2^{m(d-L)}}{\sqrt{\WW(4^m;\xi)}(1+2^m\|x-\xi\|)^{\sigma-2|\alpha|-2d}}.
\end{align*}
Since, in general, $\mu(R_\xi) \le c2^{-md/3}\WW(4^m;\xi)$
and $L$ can be arbitrarily large the above again leads to (\ref{Phi*psi}).
$\qed$

\medskip


\noindent
{\bf Proof of Lemma \ref{lem:sum<M}.}
Denote
\begin{equation}\label{def-b-diam}
h_j^\dm(x):= \sum_{\eta\in\cX_j}\frac{|h_\eta|}{(1+2^jd(x, R_\eta))^\kappa},
\quad \kappa:=\lambda-(2|\alpha|+d)|\rho|/d,
\end{equation}
where $d(x, E):=\inf_{y\in E}\|x-y\|$ is the $\ell^\infty$ distance of $x$ from $E$.
We will show that
\begin{equation}\label{bj-dm<M}
h_j^\dm(x) \le c \cM_t\Big(\sum_{\omega\in\cX_j}|h_\omega|\ONE_{R_\omega}\Big)(x), \quad x\in \RR^d_+.
\end{equation}
Evidently,
$h_j^*(x)\le h_j^\dm(x)$, $x\in \RR^d_+$,
and hence (\ref{bj-dm<M}) implies (\ref{s*<M}).
On the other hand, using (\ref{W<W}) we have for $\xi\in\cX_j$
$$
\WW(4^j;\xi)^{-\rho/d}h_\xi^*
\le\sum_{\eta\in\cX_j}\frac{\WW(4^j;\eta)^{-\rho/d}|h_\eta|}
{(1+2^j\|\xi-\eta\|)^{\lambda-(2|\alpha|+d)|\rho|/d}}
\le cH_j^\dm(x)
\quad\mbox{for $x\in R_\xi$}
$$
where $H_\eta:=\WW(4^j;\eta)^{-\rho/d}h_\eta$.
Therefore, (\ref{bj-dm<M}) yields (\ref{sum<M}) as well.

By the definition of $Q_j$ in (\ref{def-Qj})
it follows that there exists a constant $\cd>0$ depending only on $d$
such that
$$
Q_j:=\cup_{\xi\in\cX_j} R_\xi \subset [0, \cd 2^j]^d.
$$

Let $x\in \R^d$. To prove (\ref{bj-dm<M}) we consider two cases for $x$.

\medskip

{\em Case 1:} $\|x\|>2\cd 2^j$.
Then $d(x, R_\eta)>\|x\|/2$ for $\eta\in\cX_j$ and hence
\begin{align}\label{est-h*}
h_j^\dm(x)
= \sum_{\eta\in\cX_j}\frac{|h_\eta|}{(1+2^jd(x, R_\eta))^\kappa}
&\le \frac{c}{(2^{j}\|x\|)^\kappa}\sum_{\eta\in\cX_j}|h_\eta| \notag\\
&\le \frac{c4^{jd\varrho}}{(2^{j}\|x\|)^\kappa}
\Big(\sum_{\eta\in\cX_j}|h_\eta|^t\Big)^{1/t},
\end{align}
where $\varrho:=1-\min\{1, 1/t\}\le 1$ and for the last estimate
we use H\"{o}lder's inequality if $t>1$ and the $t$-triangle inequality
if $t < 1$.

Denote $Q_x:=[0, \|x\|]^d$. Evidently, $\mu(Q_x)\sim \|x\|^{2(|\alpha|+d)}$
and combining this with (\ref{est-h*}) we arrive at
\begin{align*}
h_j^\dm(x)
&\le \frac{c4^{jd}\|x\|^{2(|\alpha|+d)/t}}{(2^{j}\|x\|)^\kappa}
\Big(\frac{1}{\mu(Q_x)}\int_{Q_x}
\Big(\sum_{\eta\in\cX_j}|h_\eta|\ONE_{R_\eta}(y)\Big)^t \w(y)dy\Big)^{1/t}\\
&\le c2^{j(2d-\kappa)}\|x\|^{2(|\alpha|+d)/t-\kappa}
\cM_t\Big(\sum_{\eta\in\cX_j}|h_\eta|\ONE_{R_\eta}\Big)(x)
\le c\cM_t\Big(\sum_{\eta\in\cX_j}|h_\eta|\ONE_{R_\eta}\Big)(x)
\end{align*}
as claimed.
Here we used the fact that $\kappa > \max\{2d, 2(|\alpha|+d)/t\}$.

\medskip

{\em Case 2:} $\|x\|\le 2\cd 2^j$.
We first subdivide the tiles
$\{R_\eta\}_{\eta\in\cX_j}$ into boxes of almost equal sides of length
$\sim 2^{-j}$.
By the construction of the tiles (see (\ref{def-Rxi})) there exists
a constant $\tilde{c}>0$ such that the minimum side of each tile $R_\eta$
is $\ge \tilde{c}2^{-j}$.
Now, evidently each tile $R_\eta$ can be subdivided into a disjoint
union of boxes $R_\theta$ with centers $\theta$ such that
$$
\theta+[-\tilde{c}2^{-j-1}, \tilde{c}2^{-j-1}]^d
\subset R_\theta \subset
\theta+[-\tilde{c}2^{-j}, \tilde{c}2^{-j}]^d.
$$
Denote by $\widetilde{\cX}_j$ the set of centers of all boxes obtained
by subdividing the tiles from $\cX_j$.
Also, set $h_\theta:=h_\eta$ if $R_\theta\subset R_\eta$.
Evidently,
\begin{equation}\label{est-b*2}
h_j^\dm(x)
:= \sum_{\eta\in\cX_j}\frac{|h_\eta|}{(1+2^jd(x, R_\eta))^\kappa}
\le \sum_{\theta\in\widetilde\cX_j}\frac{|h_\theta|}{(1+2^jd(x, R_\theta))^\kappa}
\end{equation}
and
\begin{equation}\label{est-h-h}
\sum_{\eta\in \cX_j}|h_\eta|\ONE_{R_\eta}
= \sum_{\eta\in \widetilde\cX_j}|h_\theta|\ONE_{R_\theta}.
\end{equation}
Denote
$Y_0:=\{\theta\in\widetilde\cX_j: 2^j\|\theta-x\|\le \tilde{c}\}$,
\begin{align*}
Y_m &:=\{\theta\in\widetilde\cX_j: \tilde{c}2^{m-1}\le 2^j\|\theta-x\|\le \tilde{c}2^m\},
\quad \mbox{and}\\
Q_m &:=\{y\in\R^d: \|y-x\|\le \tilde{c}(2^m+1)2^{-j}\}, \quad m\ge 1.
\end{align*}
Clearly,
$\#Y_m \le c2^{md}$,
$\cup_{\theta\in Y_m}R_\theta \subset Q_m$, and
$\widetilde\cX=\cup_{m\ge 0} Y_m$.
Similarly as in (\ref{est-h*})
\begin{align*}
&\sum_{\theta\in Y_m}\frac{|h_\theta|}{(1+2^jd(x, R_\theta))^\kappa}
\le c2^{-m\kappa}\sum_{\theta\in Y_m}|h_\theta|
\le c2^{-m\kappa}2^{md\varrho}
\Big(\sum_{\theta\in Y_m}|h_\theta|^t\Big)^{1/t}\\
&\quad\le c2^{-m(\kappa-d)}
\Big(\int_{\cup_{\theta\in Y_m}R_\theta}
\sum_{\theta\in Y_m} \mu(R_\theta)^{-1}|h_\theta|^t\ONE_{R_\theta}(y)\w(y)dy\Big)^{1/t}\\
&\quad\le c2^{-m(\kappa-d)}
\Big(\frac{1}{\mu(Q_m)}\int_{Q_m}
\Big(\sum_{\theta\in Y_m}\Big(\frac{\mu(Q_m)}{\mu(R_\theta)}\Big)^{1/t}
|h_\theta|\ONE_{R_\theta}(y)\Big)^t \w(y)dy\Big)^{1/t}.
\end{align*}
Using (\ref{size-muQ}) and that $\cup_{\theta\in Y_m}R_\theta\subset Q_m$ we get
\begin{align*}
\frac{\mu(Q_m)}{\mu(R_\theta)}
&\le c\frac{2^{(m-j)d}}{2^{-jd}}
\prod_{l=1}^d \left(\frac{x_l+2^{m-j}}{\theta_l+2^{-j}}\right)^{2\alpha_j+1}\\
&\le c2^{md}\prod_{l=1}^d
\left(\frac{\theta_l+2\cdot 2^{m-j}}{\theta_l+2^{-j}}\right)^{2\alpha_j+1}
\le c2^{m(2|\alpha|+3d)}.
\end{align*}
Therefore,
$$
\sum_{\theta\in Y_m}\frac{|h_\theta|}{(1+2^jd(x, R_\theta))^\kappa}
\le c2^{-m(\kappa-d-(2|\alpha|+3d)/t)}
\cM_t\Big(\sum_{\eta\in \cX_j}|h_\eta|\ONE_{R_\eta}\Big)(x).
$$
Summing up over $m\ge 0$, taking into account that
$\kappa > d+(2|\alpha|+3d)/t$, and also using (\ref{est-b*2})
we arrive at (\ref{bj-dm<M}).
$\qed$


\medskip

\noindent
{\bf Proof of Lemma \ref{lem:a*=b*}.}
For this proof we will need an additional lemma.


\begin{lem}\label{lem:g-g}
Let $g \in \V_{4^j}$.
For any $\sigma>0$ and $L>0$
we have for $x', x''\in 2R_\xi$, where $\xi \in\cX_j$, $j\ge 0$,
\begin{equation}\label{g-g1}
|g(x')-g(x'')|\le c2^j|x'-x''|\sum_{\eta \in \cX_j}
\frac{|g(\eta)|}{(1+2^j\|\xi-\eta\|)^\sigma}
\end{equation}
and
\begin{equation}\label{g-g2}
|g(x')-g(x'')|\le c^*2^{-jL}|x'-x''|\sum_{\eta \in \cX_j}
\frac{|g(\eta)|}{(1+2^j\|\xi-\eta\|)^\sigma},
\;\; \mbox{if}\;\; \|\xi\| > (1+2\delta)\sqrt{6}\cdot 2^j.
\end{equation}
Here
$c$ and $c^*$ depend on  $\alpha$, $d$, $\delta$, and $\sigma$
and $c^*$ depends on $L$ as well;
$2R_\xi\subset \RR^d$ is the set obtained by dilating $R_\xi$
by a factor of 2 and with the same center.
\end{lem}

\noindent
{\bf Proof.}
Let $\Lambda_{4^j}$ be the kernel from (\ref{def.Ln})
with $n=4^j$, where $\ha$ is admissible of type (a)
with $v:= \delta$.
Then
$\Lambda_{4^j}*g=g$
and  $\Lambda_{4^j}(x, \cdot)\in V_{[(1+\delta)4^j]}$.
Note that $[(1+\delta)4^j]+4^j \le 2n_j-1$.
Therefore, by Corollary~\ref{cor:cubature}
$$
g(x)=\int_{\R^d} \Lambda_{4^j}(x, y)g(y)\w(y)dy
= \sum_{\eta\in\cX_j}c_\eta\Lambda_{4^j}(x, \eta)g(\eta),
$$
where $c_\eta\sim |R_\eta|\WW(4^j;\eta)$ .
From this, we have for $x', x''\in 2R_\xi$, $\xi\in\cX_j$,
\begin{align}\label{est-g-g}
|g(x')-g(x'')|
&\le \sum_{\eta\in\cX_j} c_\eta
|\Lambda_{4^j}(x', \eta)-\Lambda_{4^j}(x'', \eta)||g(\eta)|\notag\\
&\le c\|x'-x''\|\sum_{\eta\in\cX_j}c_\eta\sup_{x\in 2R_\xi}
\|\nabla \Lambda_{4^j}(x, \eta)\||g(\eta)|.
\end{align}
Note that
$(6(1+\delta)4^j+3\alpha+3)^{1/2} \le (1+\delta)\sqrt{6}\cdot 2^{j}$
for sufficiently large $j$ (depending on $\alpha$ and $\delta$).
Therefore, using Theorem~\ref{thm:est-Lambda-n} we have for $\eta\in\cX_j$
\begin{equation}\label{est-nabla1}
\|\nabla \Lambda_{4^j}(x, \eta)\|
\le \frac{c2^{j(d+1)}}{\sqrt{\WW(4^j; x)}\sqrt{\WW(4^j; \eta)}(1+2^j\|x-\eta\|)^\sigma},
\quad x\in\RR^d_+,
\end{equation}
and for any $L>0$
\begin{equation}\label{est-nabla2}
\|\nabla \Lambda_{4^j}(x, \eta)\|
\le \frac{c2^{-jL}}{(1+2^j\|x-\eta\|)^\sigma},
\quad\mbox{if} \;\; \min\{\|x\|, \|\eta\|\} > (1+\delta)\sqrt{6}\cdot 2^j.
\end{equation}

Suppose $\|\xi\|\le (1+2\delta)\sqrt{6}\cdot 2^j$ and denote
$\cX_j':= \{\eta\in\cX_j: \|\eta\|\le (1+\delta)\sqrt{6}\cdot 2^j\}$
and $\cX'':= \cX_j\setminus \cX_j'$.
We split the sum in (\ref{est-g-g}) over $\cX'$ and $\cX''$ to obtain
\begin{align*}
|g(x')-g(x'')|\le c\|x'-x''\|\Big(\sum_{\eta\in\cX_j'}\dots + \sum_{\eta\in\cX_j''}\dots\Big)
=: c\|x'-x''\|(\Sigma_1+\Sigma_2).
\end{align*}
Using (\ref{est-nabla1}), (\ref{W<W}), and that
$c_\eta\sim 2^{-jd}\WW(4^j;\eta)$ for $\eta\in\cX_j'$,
we get
\begin{align}\label{estimate-Sigma1}
\Sigma_1
&\le c2^j\sum_{\eta\in\cX_j'}\sup_{x\in 2R_\xi}
\Big(\frac{\WW(4^j; \eta)}{\WW(4^j; x)}\Big)^{1/2}
\frac{|g(\eta)|}{(1+2^j\|x-\eta\|)^\sigma}\\
&\le c2^j\sum_{\eta\in\cX_j'}
\frac{|g(\eta)|}{(1+2^j\|\xi-\eta\|)^{\sigma-2(|\alpha|+d)}}\notag
\end{align}
To estimate $\Sigma_2$ we use (\ref{est-nabla2}) and the rough estimate
$c_\eta\le c2^{jd}$.
We get
\begin{align}\label{estimate-Sigma2}
\Sigma_2
\le c2^{-j(L-d-2\sigma/3)}\sum_{\eta\in\cX_j''}
\frac{|g(\eta)|}{(1+2^j\|\xi-\eta\|)^{\sigma}}
\end{align}
Here we also used that
$$
1+2^j\|\xi-\eta\|\le 1+2^j(c2^{-j/3}+\|x-\eta\|)
\le c2^{2j/3}(1+2^j\|x-\eta\|)
\quad \mbox{for $x\in2R_\xi$}.
$$
Estimates (\ref{estimate-Sigma1}) (with sufficiently large $\sigma$)
and (\ref{estimate-Sigma2}) (with $L\ge d+2\sigma/3$)
imply (\ref{g-g1}).

In the case $\|\xi\|> (1+2\delta)\sqrt{6}\cdot 2^j$, we have
$2R_\xi \subset \{x\in \RR^d_+: \|x\|\ge (1+\delta)\sqrt{6}\cdot 2^j\}$
for sufficiently large $j$
and one proceeds just as above but uses only (\ref{est-nabla2}) as in the estimation of $\Sigma_2$.
We skip the details.
$\qed$


\medskip

We now proceed with the prove Lemma~\ref{lem:a*=b*}.
Let $g\in V_{4^j}$.
Let $\ell \ge 1$ be sufficiently large (to be determined later on) and denote
for $\xi\in\cX_j$
\begin{equation}\label{def-Xj-xi}
\cX_{j+\ell}(\xi):=\{\eta\in\cX_{j+\ell}: R_\eta\cap R_\xi \ne \emptyset\}
\quad \mbox{and}
\end{equation}
\begin{equation}\label{def-d-xi}
d_\xi:=\sup\{|g(x')-g(x'')|:
x',x''\in R_\eta \;\; \mbox{for some} \;\;\eta\in\cX_{j+\ell}(\xi)\}.
\end{equation}

Our first step is to estimate $d_\xi$, $\xi\in\cX_j$. Two cases are to be considered here.


{\em Case I:} $\|\xi\|\le (1+3\delta)\sqrt{6}\cdot 2^j$.
By (\ref{size-R-xi1})
\begin{equation}\label{incl-R}
R_\xi\sim \xi+[-2^{-j}, 2^{-j}]^d
\;\;\mbox{and}\;\;
R_\eta\sim \eta+[-2^{-j-\ell}, 2^{-j-\ell}]^d,
\;\eta\in \cX_{j+\ell}(\xi).
\end{equation}
Hence, for sufficiently large $\ell$
($\ell=\ell(d, \delta)$) we have
$\cup_{\eta\in\cX_{j+\ell}(\xi)} R_\eta \subset 2R_\xi$.
Now, using estimate (\ref{g-g1}) of Lemma~\ref{lem:g-g} with $\sigma\ge \lambda$
and the fact that $\diam (R_\eta) \sim 2^{-j-\ell}$ for $\eta\in\cX_{j+\ell}(\xi)$,
we get
\begin{equation}\label{est-d-xi1}
d_\xi \le c2^{-\ell}\sum_{\eta\in\cX_j}\frac{|g(\eta)|}{(1+2^j\|\xi-\eta\|)^\lambda},
\end{equation}
where $c>0$ is a constant independent of $\ell$.


{\em Case II:} $\|\xi\|>(1+3\delta)\sqrt{6}\cdot 2^j$.
By (\ref{size-R-xi1}) it follows that
$\|x\| > (1+2\delta)\sqrt{6}\cdot 2^j$ for
$x\in \cup_{\eta\in\cX_{j+\ell}(\xi)} R_\eta$
if $j$ is sufficiently large.
We apply estimate (\ref{g-g2}) of Lemma~\ref{lem:g-g} with $\sigma\ge \lambda$ and
$L=1$ to obtain
\begin{equation}\label{est-d-xi2}
d_\xi \le c2^{-j}\sum_{\eta\in\cX_j}\frac{|g(\eta)|}{(1+2^j\|\xi-\eta\|)^\lambda}.
\end{equation}

We next estimate $M_\xi^*$, $\xi\in\cX_j$ (see (\ref{def-s})).
Two cases for $\xi$ occur here.

\smallskip

{\em Case 1:} $\|\xi\|\le (1+4\delta)\sqrt{6}\cdot 2^j$.
Note that (\ref{incl-R}) is again valid.
By the definition of $d_\xi$ in (\ref{def-d-xi}) it follows that
$M_\xi \le m_\omega+d_\xi$ for some $\omega\in\cX_{j+\ell}(\xi)$
and hence, using (\ref{incl-R}),
$$
M_\xi \le c\sum_{\omega\in\cX_{j+\ell}}\frac{m_\omega}{(1+2^{j+\ell}\|\xi-\omega\|)^\lambda}
+ d_\xi =: \tm_\xi+ d_\xi,
\quad c=c(d, \delta, \lambda, \ell).
$$
Consequently,
\begin{equation}\label{Est-M*}
M_\xi^* \le \tm_\xi^*+ d_\xi^*.
\end{equation}
Denote $\cX_j':= \{\eta\in\cX_j: \|\eta\|\le (1+3\delta)\sqrt{6}\cdot 2^j\}$
and
$\cX_j'':= \cX_j\setminus \cX_j'$.
Now, we use (\ref{est-d-xi1})-(\ref{est-d-xi2}) to obtain
\begin{align*}
d_\xi^*
:=\sum_{\eta\in\cX_j}\frac{d_\eta}{(1+2^j\|\xi-\eta\|)^\lambda}
&\le c2^{-\ell}\sum_{\eta\in\cX_j}\sum_{\omega\in\cX_j'}
\frac{|g(\omega)|}{(1+2^j\|\xi-\eta\|)^\lambda(1+2^j\|\eta-\omega\|)^\lambda}\\
&+c2^{-j}\sum_{\eta\in\cX_j}\sum_{\omega\in\cX_j''}
\frac{|g(\omega)|}{(1+2^j \|\xi-\eta\|)^\lambda (1+2^j\|\eta-\omega\|)^\lambda}.
\end{align*}
Replacing $\cX_j'$ and $\cX_j''$ by $\cX_j$ above and shifting the order of summation
we get
\begin{align}\label{Est-d-xi}
d_\xi^*
&\le c(2^{-\ell}+2^{-j})\sum_{\omega\in\cX_j}|g(\omega)| \sum_{\eta\in\cX_j}
\frac{1}{(1+2^j\|\xi-\eta\|)^\lambda(1+2^j\|\eta-\omega\|)^\lambda}\\
&\le c(2^{-\ell}+2^{-j})\sum_{\omega\in\cX_j}
\frac{|g(\omega)| }{(1+2^j\|\xi-\omega\|)^\lambda}
\le c(2^{-\ell}+2^{-j})M_\xi^*.\notag
\end{align}
Here the constant $c$ is independent of $\ell$ and $j$,
and we used that
\begin{align}
\sum_{\eta\in\cX_j}\frac{1}{(1+2^j\|\xi-\eta\|)^\lambda(1+2^j\|\eta-\omega\|)^\lambda}
&\le \frac{c}{(1+2^j\|\xi-\omega\|)^\lambda}
\quad (\lambda >d). \label{sum<int}
\end{align}
This estimate easily follows from the fact that
$\|\xi'-\xi''\|\ge c2^{-j}$ for all $\xi', \xi''\in \cX_j$.

To estimate $\tm_\xi^*$ we use again (\ref{size-R-xi1}) and (\ref{sum<int}). We get
\begin{align*}
\tm_\xi^*
&:=\sum_{\eta\in\cX_j}\frac{\tm_\eta}{(1+2^j\|\xi-\eta\|)^\lambda}
\le c\sum_{\eta\in\cX_j}\sum_{\omega\in\cX_{j+\ell}}
\frac{m_\omega}{(1+2^j\|\xi-\eta\|)^\lambda(1+2^j\|\eta-\omega\|)^\lambda}\\
&\le c\sum_{\omega\in\cX_{j+\ell}}m_\omega\sum_{\eta\in\cX_j}
\frac{1}{(1+2^j \|\xi-\eta\|)^\lambda(1+2^j\|\eta-\omega\|)^\lambda}\\
&\le c \sum_{\omega\in\cX_{j+\ell}}
\frac{m_\omega}{(1+2^j \|\xi-\omega\|)^\lambda}
\le c2^{\ell\lambda} \sum_{\omega\in\cX_{j+\ell}}
\frac{m_\omega}{(1+2^{j+\ell}\|\theta-\omega\|)^\lambda}
=cm_\theta^*
\end{align*}
for each $\theta \in \cX_{j+\ell}(\xi)$.
Combining this with (\ref{Est-M*})-(\ref{Est-d-xi}) we obtain
$$
M_\xi^*\le c_1 m_\theta^*+ c_2(2^{-\ell}+2^{-j})M_\xi^*
\quad\mbox{for $\theta \in \cX_{j+\ell}(\xi)$,}
$$
where $c_2>0$ is independent of $\ell$ and $j$.
Choosing $\ell$ and $j$ sufficiently large (depending only on $d$, $\delta$, and $\lambda$)
this yields $M_\xi^*\le cm_\theta^*$ for all $\theta \in \cX_{j+\ell}(\xi)$.
For $j\le c$ this relation follows as above but using only (\ref{g-g1})
and taking $\ell$ large enough. We skip the details.
Thus we have shown (\ref{a*=b*}) in Case 1.

\smallskip


{\em Case 2:}  $\|\xi\| > (1+4\delta)\sqrt{6}\cdot 2^j$.
Choose $\ell\ge 1$ the same as in Case 1.
Clearly, for sufficiently large $j$ (depending only on $d$ and $\delta$)
$
\|x\| \ge (1+3\delta)\sqrt{6}\cdot 2^j
$
for $x\in\cup_{\eta\in\cX_{j+\ell}(\xi)} R_\eta$.
Hence, using (\ref{g-g2}) with $L=1$, we have
$$
M_\xi\le m_\omega+c2^{-j}\sum_{\eta\in\cX_j}\frac{|g(\eta)|}{(1+2^j\|\xi-\eta\|)^\lambda}
\le m_\omega +c2^{-j}M_\xi^*
\quad\mbox{for all $\omega\in\cX_{j+\ell}(\xi)$,}
$$
where $c>0$ is independent of $j$.
Fix $\theta\in\cX_{j+\ell}(\xi)$ and for each $\eta\in\cX_j$, $\eta\ne \xi$,
choose $\omega_\eta\in \cX_{j+\ell}(\eta)$ so that
$
\|\theta-\omega_\eta\|=\min_{\omega\in\cX_{j+\ell}(\eta)} \|\theta-\omega\|.
$
Then from above
\begin{equation}\label{Est-M*2}
M_\xi^* \le \sum_{\eta\in\cX_j}\frac{m_{\omega_\eta}}{(1+2^j\|\xi-\eta\|)^\lambda}
+c2^{-j}\sum_{\eta\in\cX_j}\frac{M_\eta^*}{(1+2^j\|\xi-\eta\|)^\lambda}
=: \Sigma_1+\Sigma_2.
\end{equation}
From (\ref{zeros-difference}) it easily follows that $\omega_\eta$ from above satisfies
$|\theta-\omega_\eta| \le c|\xi-\eta|$
and hence
\begin{equation}\label{Est-A1}
\Sigma_1\le c\sum_{\eta\in\cX_j}\frac{m_{\omega_\eta}}{(1+2^j\|\theta-\omega_\eta\|)^\lambda}
\le c2^{\ell\lambda}\sum_{\omega\in\cX_{j+\ell}}\frac{m_{\omega}}{(1+2^{j+\ell}\|\theta-\omega\|)^\lambda}
\le c_1m_\theta^*.
\end{equation}
On the other hand, using Definition~\ref{def-s*} and (\ref{sum<int}), we have
\begin{align*}
\Sigma_2
&\le c2^{-j}\sum_{\eta\in\cX_j}\sum_{\omega\in\cX_j}
\frac{M_\omega}{(1+2^j\|\xi-\eta\|)^\lambda(1+2^j\|\eta-\omega\|)^\lambda}\\
&\le c2^{-j}\sum_{\omega\in\cX_j} M_\omega\sum_{\eta\in\cX_j}
\frac{1}{(1+2^j\|\xi-\eta\|)^\lambda (1+2^j\|\eta-\omega\|)^\lambda}\\
&\le c_22^{-j}\sum_{\omega\in\cX_j}
\frac{M_\omega}{(1+2^j\|\xi-\omega\|)^\lambda}
= c_22^{-j}M_\omega^*
\end{align*}
with $c_2>0$ independent of $j$.
Combining this with (\ref{Est-M*2})-(\ref{Est-A1}) we arrive at
$$
M_\xi^*\le c_1m_\theta^*+ c_22^{-j}M_\xi^*
\quad\mbox{for $\theta\in\cX_{j+\ell}(\xi)$.}
$$
Choosing $j$ sufficiently large we get
$
M_\xi^*\le c_1m_\theta^*
$
for each $\theta\in\cX_{j+\ell}(\xi)$.
For $j\le c$ this estimate follows as in Case~1 but using only (\ref{g-g1}).
This completes the proof of Lemma~\ref{lem:a*=b*}.
\qed

\medskip


\noindent
{\bf Proof of Lemma~\ref{l:half_shannon}.}
Let $g\in \V_{4^j}$ and $0<p<\infty$.
We will utilize Definition~\ref{def-s*} and
Lemmas~\ref{lem:sum<M}-\ref{lem:a*=b*}.
To this end we select $0< t <p$ and $\lambda$ as in Definition~\ref{def-s*}.
Set $M_\xi:=\sup_{x\in R_\xi}|g(x)|$, $\xi\in\cX_j$,
and $m_\eta:=\inf_{x\in R_\eta}|g(x)|$, $\eta\in\cX_{j+\ell}$,
where $\ell\ge 1$ is the constant from Lemma \ref{lem:a*=b*}.
By (\ref{def-Wn}) and the properties of the tiles $R_\xi$
from (\ref{size-R-xi1})-(\ref{size-R-xi3}) it readily follows that
$\WW(4^{j+\ell}; y)\sim \WW(4^j,\xi)$ for $y\in R_\xi$.
We now use this, Lemmas~\ref{lem:sum<M}-\ref{lem:a*=b*} and the maximal inequality (\ref{max-ineq})
to obtain
\begin{align*}
&\Big(\sum_{\xi\in \cX_j}\WW(4^j; \xi)^{-\rho p/d}\max_{x\in R_\xi}|g(x)|^p \mu(R_\xi)\Big)^{1/p}
\le\Big\|\sum_{\xi\in \cX_j}\WW(4^j; \xi)^{-\rho/d}M_\xi^* \ONE_{R_\xi}\Big\|_p\\
&\le c\Big\|\sum_{\eta\in \cX_{j+\ell}}\WW(4^{j+\ell}; \eta)^{-\rho/d}m_\eta^* \ONE_{R_\eta}\Big\|_p
\le c\Big\|\cM_t\Big(\sum_{\eta\in \cX_{j+\ell}}\WW(4^{j+\ell}; \eta)^{-\rho/d}m_\eta \ONE_{R_\eta}\Big)\Big\|_p\\
&\le c\Big\|\sum_{\eta\in \cX_{j+\ell}}\WW(4^{j+\ell}; \eta)^{-\rho/d}m_\eta \ONE_{R_\eta}\Big\|_p
\le c\|\WW(4^j; \cdot)^{-\rho/d} g(\cdot)\|_p.
\qed
\end{align*}

\end{document}